\documentclass[final,onefignum,onetabnum]{siamonline220329}
\usepackage[mathscr]{euscript}
 
\usepackage{enumitem}
\usepackage{pgfplots}
\pgfplotsset{compat=newest}
\usetikzlibrary{plotmarks,quotes}
\usetikzlibrary{arrows.meta}
\usepgfplotslibrary{patchplots}
\newtheorem{example}{\bf{Example}}

\usepackage{lipsum}
\usepackage{amsfonts}
\usepackage{graphicx}
\usepackage{epstopdf}
\usepackage{algorithmic}
\ifpdf
  \DeclareGraphicsExtensions{.eps,.pdf,.png,.jpg}
\else
  \DeclareGraphicsExtensions{.eps}
\fi

\usepackage{enumitem}
\setlist[enumerate]{leftmargin=.5in}
\setlist[itemize]{leftmargin=.5in}

\newsiamremark{remark}{Remark}
\newsiamremark{hypothesis}{Hypothesis}
\crefname{hypothesis}{Hypothesis}{Hypotheses}
\newsiamthm{claim}{Claim}

\headers{Accuracy and componentwise accuracy in multilinear PageRank}{M. Najafi Kalyani and F. Poloni}

\title{Accuracy and componentwise accuracy in multilinear PageRank\thanks{Submitted to the editors DATE.
\funding{This research has been supported by ICSC--Centro Nazionale di Ricerca in High Performance Computing, Big Data, and Quantum Computing funded by European Union--NextGenerationEU. FP is a member of INDAM (Istituto Nazionale di Alta Matematica).}}}

\author{Mehdi Najafi Kalyani\thanks{Department of Computer Science, University of Pisa, Italy
  (\email{mehdi.najafi@di.unipi.it}, \email{federico.poloni@unipi.it}).}
\and Federico Poloni\footnotemark[2]}

\usepackage{amsopn}

\usepackage{amsmath}
\usepackage{xcolor}
\definecolor{mycolor1}{rgb}{1.00000,0.00000,1.00000}%
\definecolor{mycolor2}{rgb}{0.00000,1.00000,1.00000}%
\definecolor{mycolor3}{rgb}{0.3,0.,0.9}%
\usepackage{mathtools,booktabs,multirow,amssymb}
\usepackage{siunitx}
\usepackage{arydshln,comment}
\DeclarePairedDelimiter{\abs}{\lvert}{\rvert}
\DeclarePairedDelimiter{\norm}{\lVert}{\rVert}
\newcommand{\geqzero}{\mathrel{\geq_0}}
\DeclareMathOperator{\adj}{adj}
\DeclareMathOperator{\offdiag}{offdiag}

\newcommand{\xx}{\mathbf{x}}
\newcommand{\yy}{\mathbf{y}}
\newcommand{\ww}{\mathbf{w}}
\newcommand{\zz}{\mathbf{z}}
\newcommand{\vv}{\mathbf{v}}
\newcommand{\hh}{\mathbf{h}}
\newcommand{\mm}{\mathbf{m}}
\renewcommand{\ss}{\mathbf{s}}
\newcommand{\PP}{\mathscr{P}}
\renewcommand{\aa}{\mathbf{a}}
\newcommand{\BB}{\mathscr{B}}
\newcommand{\ec}{e_{\mathrm{cw}}}
\newcommand{\en}{e_{\mathrm{norm}}}

\ifpdf
\hypersetup{
  pdftitle={Accuracy and componentwise accuracy in multilinear PageRank},
  pdfauthor={M. Najafi Kalyani and F. Poloni}
}
\fi

\begin{document}

\maketitle

\begin{abstract}
    We study the stability with respect to perturbations and the accuracy of numerical algorithms for computing solutions to the multilinear PageRank problem $\xx = (1-\alpha)\vv + \alpha \PP \xx^2$. Our results reveal that the solution can be more stable with respect to perturbations and numerical errors with respect to the classical bounds for nonlinear systems of equations (based on the norm of the Jacobian). In detail, one can obtain bounds for the minimal solution which ignore the singularity of the problem for $\alpha=1/2$, and one can show that the limiting accuracy of the Newton method depends not on the norm of the Jacobian but on a quantity that can be much smaller thanks to the nonnegativity structure of the equation. For the minimal solution, we also suggest subtraction-free modifications to the existing algorithms to achieve componentwise stability.\\
    Some of the theoretical results we obtain are interesting even outside the scope of this problem: bounds for more general quadratic vector equations, and a partial inverse for M-matrices which remains bounded when the matrix to invert approaches singularity.
\end{abstract}



\section{Introduction}
In this work, we consider the equation
	\begin{equation}\label{eq01}
\xx =\aa+\BB\xx^2, 
\end{equation}
where $\aa \geq \mathbf{0}$ is an $n$-dimensional vector, and $\BB \geq 0$ is an $n\times n\times n$ tensor, and
\[
{(\BB\xx^2)}_i=\sum_{j,k=1}^{n}b_{ijk}x_{j}x_{k},\quad i=1,\dots,n.
\]
This equation has been studied in detail in~\cite{qve}, and it arises in the study of higher-order Markov chains; one of the applications is the so-called \emph{multilinear PageRank}~\cite{gleich2015}, in which $\aa = (1-\alpha)\vv$, $\BB = \alpha \PP$, $\alpha \in (0,1)$,  i.e.,
\begin{equation}\label{mlpr}
	{\xx} =(1-\alpha){\vv}+\alpha\PP{{\xx}}^2, 
\end{equation}
where $\vv$ and $\xx$ are stochastic (i.e., $\mathbf{1}_n^T \vv = 1$, $\mathbf{1}_n^T \xx = 1$), and $\PP$ is an $n$ dimensional stochastic tensor (or called transition probability tensor), i.e.,
\[
{p}_{ijk}\geq 0,\qquad\sum_{i=1}^{n}{p}_{ijk}=1.
\]
The theoretical properties of~\eqref{mlpr} and some applications are explored in~\cite{benson, benson2, gleich2015}. Several effective algorithms have been proposed for its numerical solution. For example, Gleich et al. \cite{gleich2015} developed Newton and fixed-point iteration methods, while Cipolla et al. \cite{cipolla2020extrapolation} used extrapolation techniques. Meini and Poloni \cite{meini2018perron} explored algorithms based on computing the Perron vector of certain nonnegative matrices, and Lai et al. \cite{lai2023anderson} applied Anderson acceleration techniques. These algorithms are currently among the most efficient for solving the equation \eqref{mlpr}.

In this paper, we study the behavior of the solutions of \eqref{eq01} and~\eqref{mlpr} under perturbations, and the accuracy at which these solutions can be computed. We are interested mostly in the multilinear PageRank problem~\eqref{mlpr}, but we formulate some of our results directly for~\eqref{eq01} for greater generality. Under two different settings, one for $\alpha \leq 1/2$ and one for $\alpha>1/2$, we show that the solution of~\eqref{mlpr} can be computed in a more stable way than one would expect based on the classical bound involving the norm of the Jacobian of the map
\begin{equation} \label{Fmap}
    F(\xx) = \aa - \BB\xx^2 - \xx.
\end{equation}

Indeed, in this problem the value $\alpha=1/2$ is an important one, since for $\alpha=1/2$ the matrix Jacobian in the solution is singular. We shall see that $\alpha=1/2$ is the switching point between two different regimes, in which the solutions have different properties.

For $\alpha \leq 1/2$, we can show stability under componentwise perturbations, i.e., using the distance defined as $d(\tilde{\xx},\xx) = \max_{i} |\xx_{i}-\tilde{\xx}_{i}| / |\xx_{i}|$. This perturbation model is natural in certain probability applications, where accuracy is required also on very small entries of the solution.

The following example shows that stability in the componentwise sense is a more stringent requirement than the usual one in the normwise sense: a perturbation can be several orders of magnitude larger when measured under the distance $d(\cdot,\cdot)$.
\begin{example}
Consider equation \eqref{eq01} with the following inputs
\[
\aa = (1-\alpha)\begin{bmatrix} 1 - \delta \\ \delta \end{bmatrix}, \quad \BB_{(1)} = \alpha 
\left[\begin{array}{cc:cc}
1 & 0.5 & 0.5 & 0 \\
0 & 0.5 & 0.5 & 1
\end{array}\right],
\]
where $\alpha \in (0,1)$ and $\delta\ll 1$ is a small positive number. It can be verified that $
\xx = \begin{bmatrix} 1 - \delta, \delta \end{bmatrix}^T$
is the solution of \eqref{eq01}.
Now, suppose $\delta$ contains a small perturbation $\varepsilon$ (e.g. due to computational error or noise in data) such that $\tilde{\delta}=\delta + \varepsilon$. We can then consider the perturbed equation
\begin{equation}\label{eq11v}
	\tilde{\xx} = \tilde{\aa}+\tilde{\BB}{\tilde{\xx}}^2
\end{equation}
with
$\tilde{\aa} = (1-\alpha)\begin{bmatrix} 1 - \tilde{\delta} ,\tilde{\delta} \end{bmatrix}^T$ and $  \tilde{\BB} = {\BB}$.
The solution of \eqref{eq11v} is $
\tilde{\xx} = \begin{bmatrix} 1 - \tilde{\delta},\tilde{\delta}\end{bmatrix}^T$.
For example, let $\delta = 10^{-6}$, $\varepsilon = 10^{-9}$. We compute the  componentwise error ($\ec$) and the normwise error ($\en$) as follows
	\[
	\ec= d(\tilde{\xx},\xx) = \max\limits_{i=1,2}\frac{|x_i - \tilde{x}_i|}{|x_i|}= \frac{\varepsilon}{\delta}=10^{-3},\qquad
	\en= \frac{{\|\xx - \tilde{\xx} \|}_2}{{\|\xx\|}_2}\approx 1.41\cdot  10^{-9}.
	\]
\end{example}

We first prove a componentwise perturbation bound, following a strategy that had been used~in \cite{xxl2} for a special case of~\eqref{eq01}. This first bound is however suboptimal, since it shows that the condition number of the problem under componentwise perturbations diverges to infinity when $\alpha\to 1/2$. We show that the bound can be improved for~\eqref{mlpr}, by restricting to perturbations that preserve the constraint $\mathbf{1}_n^T\vv = 1$ and $\mathbf{1}_n^T{\PP}_{(1)}=\mathbf{1}_{n^2}^T$. This strategy produces a bound that does not diverge and reflects the numerical behavior of the solution. 

We follow up these theoretical bounds with an algorithm to compute the minimal solution of~\eqref{mlpr} in a componentwise accurate fashion. The algorithm falls under the framework of \emph{subtraction-free algorithms}, which ensure accuracy by avoiding subtraction of quantities with the same sign~\cite{alfa,Ocinneide,xxl2,xue2012accurate}.
O'Cinneide~\cite{Ocinneide} and Alfa et al. \cite{alfa} demonstrated that the inverse of a nonsingular M-matrix can be determined accurately using a triplet representation and the \emph{GTH algorithm}, a variant of Gaussian elimination. Our approach combines the GTH algorithm with an accurately crafted modification of either the Newton method or the block Jacobi method, to obtain a numerical solution for \eqref{mlpr} with high componentwise accuracy.



When $\alpha > 1/2$, unfortunately the componentwise stability of a stochastic solution $\ss$ cannot be guaranteed: indeed, examples are known where a stochastic solution varies in a non-Lipschitz way or disappears under perturbations, see~\cite{meini2018perron}. However, we show that the normwise error of the solution computed by the Newton method can be several orders of magnitude smaller than what is predicted by the standard bounds for nonlinear equations based on $\norm{R_{\ss}^{-1}}$ (see e.g.~\cite{tisseur} for these bounds).

We present numerical experiments to demonstrate the tightness of our bounds and the efficiency and reliability of our algorithms.

The paper is organized as follows. In~\cref{sec:results}, we summarize some notations, definitions and known results that we will need in the rest of the paper. In \cref{sec:grahM}, we present a subtraction-free representation for the inverses of M-matrices. Our main theoretical bounds concerning the componentwise error are detailed in \cref{sec:main}. In \cref{sec:algorithms}, we suggest various algorithms to compute the minimal solution $\mm$ with componentwise accuracy. In~\cref{sec:limitaccuracy}, we study instead the limiting accuracy of Newton's method for stochastic solutions. These results are followed by a discussion of experimental results in \cref{sec:experiments}. Finally, we add some concluding remarks in \cref{sec:conclusions}.
\section{Useful results} \label{sec:results}

\subsection{Tensor notation} Let $\mathscr{A}$ be an $n\times n\times n$ tensor. The tensor-vector product between $\mathscr{A}$ and the vectors $\xx$ and $\yy$ is defined by
\[
{(\mathscr{A}\xx\yy)}_i=\sum_{j,k=1}^{n}a_{ijk}x_{j}y_{k},\quad i=1,\dots,n.
\]
The product in equation \eqref{eq01} is a special case of the above product where we set $\yy=\xx$. We consider the matrices associated with the linear maps
$\yy\to \mathscr{A}\xx\yy$ and $\yy\to \mathscr{A}\yy\xx$,
which are denoted as
$\mathscr{A}{\xx{:}}$ and $\mathscr{A}{{:}\xx}$, respectively.
These matrices are given by
\[
{(\mathscr{A}{\xx{:}})}_{ij}=\sum_{k=1}^{n}a_{ikj}x_{k},\qquad
{(\mathscr{A}{{:}\xx})}_{ij}=\sum_{k=1}^{n}a_{ijk}x_{k}.
\]

In particular, we will use the matrix
\begin{equation} \label{Rx}
R_{\xx} = I_n - \BB\xx{:} - \BB{:}\xx,    
\end{equation}
i.e., the unique matrix such that $R_{\xx}\hh = \hh - \BB\xx\hh - \BB\hh\xx$ for all vectors $\hh$. This matrix is the negative of the Jacobian of the nonlinear map~\eqref{Fmap}.

For the multilinear PageRank~\eqref{mlpr}, it holds that
\begin{align*}
\mathbf{1}_n^T R_{\mathbf{\xx}} &= \mathbf{1}_n^T - \mathbf{1}_n^T\BB\xx{:} - \mathbf{1}_n^T\BB{:}\xx = \mathbf{1}_n^T - \alpha (\mathbf{1}_n^T\xx) \mathbf{1}_n^T - \alpha (\mathbf{1}_n^T\xx) \mathbf{1}_n^T = (1-2\alpha (\mathbf{1}_n^T\xx)) \mathbf{1}_n^T.
\end{align*}
In particular, when $\xx$ is a stochastic solution of~\eqref{mlpr}, we have $\mathbf{1}_n^T R_{\mathbf{\xx}} = (1-2\alpha)\mathbf{1}_n^T$, which shows that $R_{\mathbf{\xx}}$ is an M-matrix if and only if $\alpha \leq 1/2$, and a nonsingular one when $\alpha < 1/2$.

For the tensor $\mathscr{A}$, the first mode unfolding (denoted as $\mathscr{A}_{(1)}$) transforms the tensor into a matrix of size $n\times n^2$, where $({\mathscr{A}_{(1)})}_{i,j+(k-1)n}=a_{ijk},\,\, i,j,k\in\{1,\dots,n\}$.

\subsection{Componentwise perturbation distance}

We are interested in the case in which $\tilde{\BB}$ is an entrywise multiplicative perturbation of $\BB$, i.e., 
\begin{equation*}
	|\tilde{b}_{ijk}-b_{ijk}|\leq \varepsilon b_{ijk},
\end{equation*}
 for all indices $i,j,k=1,2,\dots,n$. Alternatively, we can define
\[
d(\tilde{\BB},\BB) := \max_{(i,j,k)} \frac{\abs{\tilde{b}_{ijk} - b_{ijk}}}{\abs{b_{ijk}}},
\]
with the convention that $0/0=0$ and $b/0=\infty$ for any $b\neq 0$, so that
\[
d(\tilde{\BB},\BB) \leq \varepsilon \iff (1-\varepsilon)\BB \leq \tilde{\BB} \leq (1+\varepsilon)\BB.
\]
In particular, this kind of perturbation model implies that the zero patterns of $\BB$ and $\tilde{\BB}$ are the same.

It is useful to note that this definition is not symmetric: $d(\tilde{\BB},\BB)\neq d(\BB,\tilde{\BB})$.  We can extend the definition of $d(\cdot,\cdot)$ to vectors and matrices in a similar way, only with a different number of indices.

\subsection{Results on quadratic vector equations}

The following theorem establishes some basic properties of~\eqref{eq01}, found in~\cite{qve}. Here and in the following, the relation symbols $\leq$ and $<$ are used in the componentwise sense for vectors and matrices.
\begin{theorem}\label{thm.21} Suppose~\eqref{eq01} has at least a  non-negative solution. Then,
\begin{enumerate}
    \item Among all nonnegative solutions of~\eqref{eq01}, there exists a \emph{minimal} one, i.e., a solution $\mm$ such that $\mm \leq \xx$ for each other solution $\xx \geq \mathbf{0}$.
    \item The matrix $R_{\mm}$, defined as in~\eqref{Rx} where $\mm$ is the minimal nonnegative solution, is an M-matrix.
    \item The fixed-point iteration 
    \begin{equation}\label{eq02}
	   {\xx}_{k+1} ={\aa}+\BB{\xx}_k^2, \qquad k = 0, 1,\dots,
    \end{equation}
    as well as the Newton method on~\eqref{Fmap}, which takes the form
    \begin{equation}\label{newton}
	   {\xx}_{k+1} = \xx_k + R_{\xx_k}^{-1}\underbrace{({\aa}+\BB{\xx}_k^2-\xx_k)}_{=F(\xx_k)} = R_{\xx_k}^{-1}({\aa}-\BB{\xx}_k^2), \qquad k = 0, 1,\dots,
    \end{equation}
    when started from $\xx_0 = \mathbf{0}$, have the properties that $F(\xx_k)\geq \mathbf{0}$ at each step, and that the iterates converge monotonically to the minimal solution $\mm$:
    \[
    \mathbf{0} = \xx_0 \leq \xx_1 \leq \dots \leq \xx_k \leq \dots \leq \mm, \quad \lim_{k\to\infty} \xx_k = \mm.
    \]
\end{enumerate}
\end{theorem}
We shall assume in our paper that $\mm > \mathbf{0}$, to exclude degenerate cases where certain entries are constantly equal to $0$ along our iterations; see also~\cite[Section~8]{qve} for more discussion of these degenerate cases. This assumption, in particular, ensures that $R_{\xx}$ is a non-singular M-matrix for all $\xx \leq \mm$, $\xx \neq \mm$.

The papers~\cite{gleich2015,meini2018perron} contain additional results for~\eqref{mlpr}, which is a special case of~\eqref{eq01}; we summarize the ones that we need in the following theorem.
\begin{theorem}
        If $\alpha \in (0,1/2]$, the minimal nonnegative solution $\mm$ of~\eqref{mlpr} is stochastic.
        
        If $\alpha \in (1/2,1)$, the minimal nonnegative solution $\mm$ of~\eqref{mlpr} is sub-stochastic, with $\mathbf{1}_n^T\mm = \frac{1-\alpha}{\alpha} < 1$. In addition to $\mm$, \eqref{mlpr}~has at least one (non-minimal) stochastic solution.
\end{theorem}

\section{A subtraction-free representation for inverses of M-matrices}\label{sec:grahM}
In this section, we explore the properties of certain M-matrices and their inverses, and define a certain partial inverse that we will need in the following.

The following result is a consequence of the \emph{all-minors matrix tree theorem} (see, e.g., \cite{Abdesselam,chaiken}).
\begin{theorem} \label{matrixtree}
    Consider a directed graph without loops on nodes $\{0,1,2,\dots,n\}$, with weights $w_{ij}$, $i\neq j$, $i,j \in \{0,1,\dots,n\}$. Let $M\in\mathbb{R}^{n\times n}$ be the submatrix obtained by removing the row and column with index $0$ from the Laplacian of this graph; i.e., $M$ is the matrix such that
    \begin{equation} \label{Mmatrixtree}
    M_{ij} = -w_{ij} \quad \text{for all  $i\neq j \in \{1,\dots,n\}$, and} \quad 
    M\mathbf{1} = \begin{bmatrix}
        w_{10}\\
        w_{20}\\
        \vdots\\
        w_{n0}
    \end{bmatrix}.
    \end{equation}
    Then,
    \begin{equation} \label{detM_smatrixtree}
    \det M = \sum_{\mathcal{T}\in \mathcal{F}} \prod_{(i,j)\in\mathcal{T}} w_{ij},        
    \end{equation}
    where $\mathcal{F}$ is a certain family of subsets of the graph edges. Moreover, for each $k,l\in\{1,2,\dots,n\}$ the $(k,l)$ entry of the adjugate matrix $\operatorname{adj} M$ is given by
    \[
    (\operatorname{adj} M)_{kl} = \sum_{\mathcal{T}\in \mathcal{F}_{kl}} \prod_{(i,j)\in\mathcal{T}} w_{ij},
    \]
    for another such family $\mathcal{F}_{kl}$.

    More precisely, $\mathcal{F}$ is the family of all spanning trees in the graph, oriented such that all edges point towards node $0$. Similarly, each element of $\mathcal{F}_{kl}$ is the union of two disjoint trees $\mathcal{T} = \mathcal{T}_0 \cup \mathcal{T}_1$, where
    \begin{enumerate}
        \item $\mathcal{T}_0 \cup \mathcal{T}_1$ together span all nodes of the graph;
        \item $\mathcal{T}_0$ contains node $0$, and all edges point towards it;
        \item $\mathcal{T}_1$ contains nodes $k,l$, and all edges point towards $l$.
    \end{enumerate}
    The family $\mathcal{F}_{kl}$ contains all such pairs of trees.
\end{theorem}
\begin{example}
    Let us consider the complete graph on $\{0,1,2,3\}$. We have
    \[
    M = \begin{bmatrix}
            w_{1,0}+w_{1,2}+w_{1,3} & -w_{1,2} & -w_{1,3}\\ -w_{2,1} & w_{2,0}+w_{2,1}+w_{2,3} & -w_{2,3}\\ -w_{3,1} & -w_{3,2} & w_{3,0}+w_{3,1}+w_{3,2} 
    \end{bmatrix}
    \]
    and
    \begin{align*}
    \det M &=     w_{1,0}\,w_{2,0}\,w_{3,0}+w_{1,0}\,w_{2,0}\,w_{3,1}+w_{1,0}\,w_{2,1}\,w_{3,0}+w_{1,0}\,w_{2,0}\,w_{3,2}+w_{1,0}\,w_{2,1}\,w_{3,1}\\ &+w_{1,2}\,w_{2,0}\,w_{3,0}+w_{1,0}\,w_{2,1}\,w_{3,2}+w_{1,0}\,w_{2,3}\,w_{3,0}+w_{1,2}\,w_{2,0}\,w_{3,1}+w_{1,3}\,w_{2,0}\,w_{3,0}\\ &+w_{1,0}\,w_{2,3}\,w_{3,1}+w_{1,2}\,w_{2,0}\,w_{3,2}+w_{1,3}\,w_{2,1}\,w_{3,0}+w_{1,2}\,w_{2,3}\,w_{3,0}+w_{1,3}\,w_{2,0}\,w_{3,2}\\ 
    &+w_{1,3}\,w_{2,3}\,w_{3,0}.        
    \end{align*}
    Each summand is the product of the weights in a spanning tree pointing towards $0$; for instance the second summand $w_{1,0}\,w_{2,0}\,w_{3,1}$ corresponds to
    
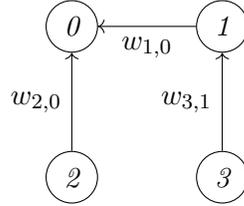
\begin{center}
    \begin{tikzpicture}
    \node[draw, circle] (A) at (0, 2) {0};
    \node[draw, circle] (B) at (2, 2) {1};
    \node[draw, circle] (C) at (0, 0) {2};
    \node[draw, circle] (D) at (2, 0) {3};

    \draw[->] (C) edge["$w_{2,0}$"] (A);
    \draw[->] (B) edge["$w_{1,0}$"] (A);
    \draw[->] (D) edge["$w_{3,1}$"] (B);

    \end{tikzpicture}
    \end{center}
    The adjugate matrix $\operatorname{adj} M$, i.e., the matrix that satisfies $M \operatorname{adj}(M) = \det(M) I_n$, has $(2,1)$ entry
    \[
    (\operatorname{adj} M)_{21} = - \det \begin{bmatrix}
        -w_{2,1} & -w_{2,3}\\ -w_{3,1} &  w_{3,0}+w_{3,1}+w_{3,2}
    \end{bmatrix} = w_{2,1}\,w_{3,0}+w_{2,1}\,w_{3,1}+w_{2,1}\,w_{3,2}+w_{2,3}\,w_{3,1}.
    \]
    Each summand is the product of the weights in a pair of trees satisfying the conditions in the theorem; for instance, the second summand $w_{2,1}\,w_{3,1}$ corresponds to
    
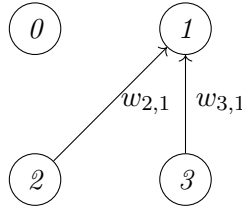
\begin{center}
    \begin{tikzpicture}
    \node[draw, circle] (A) at (0, 2) {0};
    \node[draw, circle] (B) at (2, 2) {1};
    \node[draw, circle] (C) at (0, 0) {2};
    \node[draw, circle] (D) at (2, 0) {3};

    \draw[->] (C) edge["$w_{2,1}$" right] (B);
    \draw[->] (D) edge["$w_{3,1}$" right] (B);
    \end{tikzpicture}
    \end{center}
    Here $\mathcal{T}_0$ is the trivial tree containing only node $0$, while $\mathcal{T}_1$ contains nodes $k=2$ and $l=1$, and its edges point towards the latter.
\end{example}
Note that this result can be applied to give an expression for the determinant and the adjoint of \emph{any} matrix $M$: indeed, given a matrix $M$ we can choose weights $w_{i,j}$ (not necessarily positive) such that~\eqref{Mmatrixtree} holds. Indeed, one sees that the opposites of the offdiagonal entries
\begin{equation} \label{offdiagonal}
    -M_{i,j} = w_{i,j}, \quad i\neq j, \quad i,j\in \{1,2,\dots,n\}    
\end{equation}
and the row sums
\begin{equation} \label{rowsums}
M\mathbf{1} = \mathbf{w} = \begin{bmatrix}
    w_{1,0}\\
    w_{2,0}\\
    \vdots\\
    w_{n,0}
\end{bmatrix}    
\end{equation}
uniquely determine the matrix $M$. (Note that the weights of the form $w_{0,j}$ are formally defined for the graph, but they do not appear in our theorem.)

In the following, we will sometimes want to emphasize the dependence on $M$ on the row sums~\eqref{rowsums}, while assuming that the offdiagonal entries~\eqref{offdiagonal} are given and fixed; in this case, we write $M = M(\mathbf{w})$.

This expression allows one to prove componentwise stability results for inverses of M-matrices. The following result appeared in~\cite{alfa,Ocinneide}, who obtained it without relying on the combinatorial interpretation that we have described.
\begin{theorem}\label{MMt}
    Let $M, \tilde{M}$ be two invertible M-matrices with $\tilde{M}\mathbf{1}\geq \mathbf{0}$, $M\mathbf{1}\geq \mathbf{0}$. Assume that
    $d(\operatorname{offdiag}(\tilde{M}), \operatorname{offdiag}(M)) \leq \varepsilon$ and $d(\tilde{M}\mathbf{1},M\mathbf{1}) \leq \varepsilon$. Then,
    \[
    d(\tilde{M}^{-1}, M^{-1}) \leq (2n-1)\varepsilon + \mathcal{O}(\varepsilon^2).
    \]
\end{theorem}
\begin{proof}
    We can associate to $\tilde{M}$ and $M$ weights $\tilde{w}_{i,j}$ and $w_{i,j}$ as in Theorem~\ref{matrixtree}; then the hypotheses imply that $w_{ij}, \tilde{w}_{ij}\geq 0$ for all $i,j$ and
    \[
    (1-\varepsilon)w_{ij} \leq \tilde{w}_{ij} \leq (1+\varepsilon)w_{ij} \quad \text{for all $i,j$}.
    \]
    Then, we have
    \[
    \det \tilde{M} = \sum_{\mathcal{T} \in\mathcal{F}} \prod_{(i,j)\in\mathcal{T}} \tilde{w}_{i,j} \leq \sum_{\mathcal{T} \in\mathcal{F}} \prod_{(i,j)\in\mathcal{T}} (1+\varepsilon)w_{i,j} = (1+\varepsilon)^n \sum_{\mathcal{T} \in\mathcal{F}} \prod_{(i,j)\in\mathcal{T}} w_{i,j} = (1+\varepsilon)^n \det M,
    \]
    since each product has $n$ nonnegative factors. Analogously one can show that $\det \tilde{M} \geq (1-\varepsilon)^n \det M$, and that
    \[
        (1-\varepsilon)^{n-1} (\adj M)_{kl} \leq 
        (\adj \tilde{M})_{kl} \leq
        (1+\varepsilon)^{n-1} (\adj M)_{kl}.
    \]
    Since $M^{-1} = \frac{1}{\det M} \adj M$, these inequalities imply
    \[
    (1 - (2n-1)\varepsilon + \mathcal{O}(\varepsilon^2))M^{-1} =
    \frac{(1-\varepsilon)^{n-1}}{(1+\varepsilon)^n} M^{-1} \leq \tilde{M}^{-1} \leq \frac{(1+\varepsilon)^{n-1}}{(1-\varepsilon)^n} M^{-1} = (1 + (2n-1)\varepsilon + \mathcal{O}(\varepsilon^2))M^{-1},
    \]
    from which the statement follows.
\end{proof}
\begin{remark}
    We believe there is a typo in~\cite[Equations~(1) and~(2)]{alfa}: the exponents $n$ and $n-1$ in the numerator and denominator are swapped. Note that the later~\cite[Theorem~2.2]{xue2012accurate} contains the correct version.
\end{remark}
With these tools, we can prove a special decomposition for the inverse of a matrix.
\begin{theorem}\label{thmRS}
    Let $M\in\mathbb{R}^{n\times n}$ be an invertible matrix with $M_{ij}  \leq 0$ for $i\neq j$, and set $\mathbf{w} = M\mathbf{1}$. There exists a decomposition $M^{-1} = R+S$ with the following properties:
    \begin{enumerate}
        \item $R$ is rank-1: $R = \mathbf{1}\mathbf{z}^T$ for a certain vector $\mathbf{z}$.
        \item $S$ is stable: if $M(\mathbf{w}_k)$ is a sequence of irreducible matrices with $\mathbf{w}_k \to \mathbf{0}$ (while the off-diagonal entries are constant), the corresponding $S(\mathbf{w}_k)$ are bounded.
        \item $S$ acts like the inverse of $M$ over a certain subspace of row vectors: if $\vv$ is a vector such that $\vv^T \mathbf{1} = 0$, then $\vv^{T} M^{-1} = \vv^T S$.
    \end{enumerate}
    Moreover, when $\mathbf{w}\geq \mathbf{0}$ (and hence $M$ is an M-matrix):
    \begin{enumerate}[resume]
        \item $R,S\geq 0$.
        \item $R,S$ are entrywise well-conditioned: given two regular M-matrices $M, \tilde{M}$ with $\tilde{M}\mathbf{1}\geq \mathbf{0}$, $M\mathbf{1}\geq \mathbf{0}$, if
    $d(\operatorname{offdiag}(\tilde{M}), \operatorname{offdiag}(M)) \leq \varepsilon$ and $d(\tilde{M}\mathbf{1},M\mathbf{1}) \leq \varepsilon$, then the decompositions $M^{-1}=R+S$ and $\tilde{M}^{-1}=\tilde{R}+\tilde{S}$ satisfy
    \begin{equation} \label{RSstability}
    d(\tilde{R}, R) \leq (2n-1)\varepsilon + \mathcal{O}(\varepsilon^2), \quad d(\tilde{S}, S) \leq (2n-1)\varepsilon + \mathcal{O}(\varepsilon^2).    
    \end{equation}
    \end{enumerate}
    \begin{proof}
        Note that we can apply Theorem~\ref{matrixtree}, and $w_{1,0},\dots,w_{n,0}$ are the entries of $\mathbf{w}$.
        Let us consider the families $\mathcal{F}_{kl}$ that appear in Theorem~\ref{matrixtree}. For each $\mathcal{F}_{kl}$, we write $\mathcal{F}_{kl} = \mathcal{R}_{kl} \cup \mathcal{S}_{kl}$, where $\mathcal{R}_{kl}$ is the family of pairs $\mathcal{T}=\mathcal{T}_0\cup \mathcal{T}+1$ in which $\mathcal{T}_0$ is the trivial tree containing only the node $0$, and $\mathcal{S}_{kl}$ is its complement. We set
        \begin{equation} \label{RS_matrixtree}
        R_{kl} = \frac{1}{\det M} \sum_{\mathcal{T}\in \mathcal{R}_{kl}} \prod_{(i,j)\in\mathcal{T}} w_{i,j}, \quad 
        S_{kl} = \frac{1}{\det M} \sum_{\mathcal{T}\in \mathcal{S}_{kl}} \prod_{(i,j)\in\mathcal{T}} w_{i,j}.            
        \end{equation}
        For all $\mathcal{T} \in \mathcal{R}_{kl}$, since $\mathcal{T}_0$ is a singleton, $\mathcal{T}_1$ contains all nodes $k\in\{1,2,\dots,n\}$; hence the same $\mathcal{T}$ belongs to $\mathcal{R}_{kl}$ for each $k\in\{1,2,\dots,n\}$. This shows that $\mathcal{R}_{1l}=\mathcal{R}_{2l}=\dots=\mathcal{R}_{nl}$ for each $l$, i.e., $R$ has the form $\mathbf{1}\mathbf{z}^T$ for a suitable nonnegative vector $\mathbf{z}$.

        We now show the boundedness of $S$. In the expansion~\eqref{detM_smatrixtree} for $\det M$, we single out the terms that contain only one factor of the form $w_{i,0}$ for some $i$, and write
        \[
            \det M = w_{1,0} D_1 + w_{2,0}D_2 + \dots + w_{n,0}D_n + \mathcal{O}(\norm{\mathbf{w}}^2),
        \]
        where $D_i$ are polynomials in the weights $w_{i,j}$ with $i,j\neq 0$, which do not depend on the entries of $\mathbf{w}$. Similarly, for a fixed $k,l$ we can find polynomials $N_i$ such that the numerator of $S_{kl}$ is
        \[
        S_{kl}\det M = \sum_{\mathcal{T}\in \mathcal{S}_{kl}} \prod_{(i,j)\in\mathcal{T}} w_{ij} = w_{1,0} N_1 + w_{2,0}N_2 + \dots + w_{n,0}N_n + \mathcal{O}(\norm{\mathbf{w}}^2).
        \]
        Let $m = m(\mathbf{w})$ be the index for which $w_{m,0}$ is largest; then,
        \begin{align*}
        \lim\sup_{\mathbf{w}\to \mathbf{0}} S_{kl} &= 
        \lim\sup_{\mathbf{w}\to \mathbf{0}} \frac{w_{1,0} N_1 + w_{2,0}N_2 + \dots + w_{n,0}N_n + \mathcal{O}(\norm{\mathbf{w}}^2)}{w_{1,0} D_1 + w_{2,0}D_2 + \dots + w_{n,0}D_n + \mathcal{O}(\norm{\mathbf{w}}^2)} \\
        & \leq \lim\sup_{\mathbf{w}\to \mathbf{0}} \frac{w_{m,0} (N_1 + N_2 + \dots + N_n) + \mathcal{O}(\norm{\mathbf{w}}^2)}{w_{m,0}D_m}\\
        & \leq \frac{N_1 + N_2 + \dots + N_n}{\min_{i=1,\dots,n} D_i}.
        \end{align*}
        To show that this expression is bounded, it is sufficient to show that $D_i> 0$ for all $i$. We are assuming that $M$ is irreducible; hence the graph on $\{1,2,\dots,n\}$ with weights $w_{i,j}$ is irreducible. In particular, thanks to irreducibility, one can find a tree (with non-zero weights) $\hat{\mathcal{T}}$ in which each node is connected to $i$ by edges pointing to $i$. Then, $\hat{\mathcal{T}} \cup (i,0)$ is one of the trees that appear in the sum $w_{i,0}D_i$. Thus this sum contains at least a positive term, and $D_i > 0$.

        If $\mathbf{w}\geq \mathbf{0}$, we have $R,S\geq 0$, as the sums in their definition~\eqref{RS_matrixtree} involve nonnegative terms: $w_{i,j}\geq 0$ for all $i,j$. 

        To prove~\eqref{RSstability}, we argue as in the proof of Theorem~\ref{MMt}; the only difference is that we have sums over $\mathcal{R}_{kl}$ and $\mathcal{S}_{kl}$ instead of the sum over $\mathcal{F}_{kl}$ that gives $\operatorname{adj}(M)$.

    \end{proof}
\end{theorem}
\section{Componentwise perturbation analysis}\label{sec:main}
The purpose of this section is to derive upper bounds for the componentwise error of equation \eqref{eq01}. 
\subsection[]{The role of $\yy$}
Let $\mm$ be the minimal solution of~\eqref{eq01}. We shall assume in the following that $R_{\mm}$ is nonsingular, which is the most common case.

Note that the minimal solution $\mm$ of~\eqref{eq01} also satisfies
\[
\mm - \BB\mm^2 - \BB\mm^2 = \aa - \BB\mm^2,
\]
i.e.,
\[
R_{\mm} \mm = \aa - \BB\mm^2
\]
or 
\begin{equation} \label{xy}
\mm = R_{\mm}^{-1}\aa - R_{\mm}^{-1}\BB\mm^2.
\end{equation}
We define
\[
\yy = R_{\mm}^{-1} \aa,
\]
i.e., $\yy$ is the unique solution of the system of linear equations
\[
\yy = \aa + \BB\mm\yy + \BB\yy\mm.
\]
Then, it follows from~\eqref{xy} that $\mm \leq \yy$, since $R_{\mm}^{-1}\BB\mm^2 \geq \mathbf{0}$.

The quantity $\yy$ is considered also in~\cite{xxl2} for a matrix equation that reduces to a special case of~\eqref{eq01}, after vectorization; there, the authors argue that $\yy$ is usually close to $\mm$, and that it is easier to bound the entries of $\hh$ in terms of those of $\yy$ rather than those of $\mm$. The following result extends one in~\cite{xxl2} to our more general case.
\begin{theorem}\label{thm0}
 The minimal solution $ {\mm}$ and
   ${\yy}$ have the same zero pattern, i.e., if $m_i=0$ then $y_i=0$ for $ i = 0, 1,\dots,n$, and viceversa.
\end{theorem}
\begin{proof}
We use the notation $\aa \geqzero \mathbf{b}$ to mean that $\aa_i = 0 \implies \mathbf{b}_i = 0$, i.e., $\mathbf{b}$ has zeros in the same positions as $\aa$, plus possibly others. Clearly $\yy \geq \mm \geq \mathbf{0}$ implies $\yy \geqzero \mm$, so it is enough to prove $\mm \geqzero \yy$.

In view of~\eqref{eq01} and nonnegativity, we have $\mm \geqzero \aa$ and $\mm \geqzero \BB\mm^2$. Consider the iteration
\[
\yy_{k+1} = \aa + \BB\mm\yy_k + \BB\yy_k\mm, \quad \yy_0 = \mathbf{0}.
\]
One can show by induction that the sequence $\yy_k$ is increasing and bounded above by $\yy$; hence we must have $\lim \yy_k = \yy$, since the limit solves the equation that defines $\yy$.

We shall show by induction on $k$ that $\mm \geqzero \yy_k$. Clearly $\mm \geqzero \yy_0 = \mathbf{0}$. Now assume that $\mm \geqzero \yy_k$. We can see that this implies $\BB\mm^2 \geqzero \BB\mm\yy_k$: indeed, 
\begin{multline}\label{mul}
(\BB{\mm}^2)_i = 0 \iff b_{ij\ell} = 0 \vee m_j=0 \vee m_\ell=0 \text{ for all $j,\ell$} \\ \implies b_{ij\ell} = 0 \vee m_j=0 \vee (\yy_k)_\ell=0 \text{ for all $j,\ell$} \iff (\BB\mm\yy_k)_i = 0.\nonumber
\end{multline}
Similarly, $\BB\mm^2 \geqzero \BB\yy_k\mm$. Thus,
\[
\mm = \aa + \BB\mm^2  \geqzero \aa + \BB\mm\yy_k + \BB\yy_k\mm = \yy_{k+1}.
\]
This completes the induction proof. Hence, also $\mm \geqzero \lim \yy_k = \yy$.
 \end{proof}

 Since $\yy$ and $\mm$ have the same zero pattern by Theorem~\ref{thm0}, we can define
\begin{equation}\label{kapavalue}
\kappa = \max_{m_i\neq 0} \frac{y_i}{m_i},
\end{equation}
so that $\yy \leq \kappa \mm$. In~\cite{xxl2}, this quantity $\kappa$ appears in the bounds. We also prove a bound for $R_{\mm}^{-1}\mm$ based on this quantity $\kappa$.
\begin{proposition}
Assume that $\yy = R_{\mm}^{-1} \aa$, where $\yy \leq \kappa \mm$, and $\kappa$ is defined as \eqref{kapavalue}. Then, we have
\begin{align*}
  R_{\mm}^{-1} \mm &\leq(2\kappa-1)\mm,  \\
  R_{\mm}^{-1}\BB\mm^2 & \leq(\kappa-1)\mm.
\end{align*}
\end{proposition}
\begin{proof}
Using equation \eqref{xy}, we can express $R_{\mm}^{-1}\mm$ as follows
\[
R_{\mm}^{-1} \mm=R_{\mm}^{-1} (\aa+{\BB}\mm^2)=\yy+R_{\mm}^{-1}\BB\mm^2=\yy+(\yy-\mm).
\]
It now follows from $\yy \leq \kappa \mm$ that
\begin{equation}\label{kapax}
\begin{split}
 R_{\mm}^{-1} \mm&=\yy+(\yy-\mm)\leq(2\kappa-1)\mm,\\
 R_{\mm}^{-1}\BB\mm^2&=\yy-R_{\mm}^{-1} \mm\leq(\kappa-1)\mm.
\end{split}
\end{equation}
\end{proof}

\subsection{Componentwise error bounds}
The following perturbation bound holds for the minimal solution of equation~\eqref{eq01} when $R_{\mm}$ is invertible, i.e., ${\rho(\BB\mm{:} + \BB{:}\mm)} < 1$, where $\rho(\cdot)$ is the spectral radius. In particular, this includes the case of~\eqref{mlpr} when $\alpha < 1/2$.
\begin{theorem}
	Let $\BB=(b_{ijk})$ and $\tilde{\BB}=(\tilde{b}_{ijk})$ be nonnegative $n$-dimensional tensors. Let $\mm$ be the minimal solution of~\eqref{eq01}, and assume that $R_{\mm}$ is invertible. Suppose that $d(\tilde{\BB},\BB) \leq \varepsilon$ 
    and $d(\tilde{\aa},{\aa}) \leq \varepsilon$, 
    for a value $\varepsilon$ small enough to ensure
    \begin{equation} \label{Mmatrixcond}
    (1+\varepsilon){\rho(\BB\mm{:} + \BB{:}\mm)} < 1.
    \end{equation}
    and
    \begin{equation}\label{condDelta}
    \varepsilon+\varepsilon^2 <\frac{1}{4\gamma^2(2\kappa-1)(\kappa-1)},    
    \end{equation}  
    where $\gamma=\frac{(1+\varepsilon)^{n-1}}{(1-\varepsilon)^{n}}$, and $\kappa$  is defined as \eqref{kapavalue}. Then, \eqref{eq11v} has a solution $\tilde{\mm}$ such that 
    \[
        d(\tilde{\mm},\mm) 
        \leq \frac{2\varepsilon(2\kappa-1)(1+\varepsilon)^{n-1}}{(1-\varepsilon)^{n}}.
	\]
\end{theorem}
\begin{proof}
    First, note that condition~\eqref{Mmatrixcond} ensures that
    \[
    \rho(\tilde{\BB}\mm{:}+\tilde{\BB}{:}\mm) \leq (1+\varepsilon)\rho(\BB\mm{:}+\BB{:}\mm) < 1,
    \]
    so $\tilde{R}_{\mm}$ is an M-matrix.

	We define $
	B_{k_{\hh}}=\left\{{\hh}\in \mathbb{R}^{n}: \abs{\hh} \leq k_{\hh} \mm\right\}$.
	For a vector $\hh \in B_{k_{\hh}}$, we have that $\tilde{\mm} = \mm+ \hh$ satisfies~\eqref{eq02} if and only if
 \begin{align*}
 \hh=\tilde{\mm}-\mm&={\tilde{\aa}-{\aa}}+\tilde{\BB}(\mm+\hh)^2-{\BB}\mm^2\\
&={\tilde{\aa}-{\aa}}+ \tilde{\BB}\mm^2+(\tilde{\BB}\mm\hh+\tilde{\BB}\hh\mm)+\tilde{\BB}\hh^2-{\BB}\mm^2.	    
 \end{align*}   
Alternatively, we get
\[
	\tilde{R}_{\mm}\hh:={\hh-\tilde{\BB}\mm\hh-\tilde{\BB}\hh\mm} =
	{\tilde{\aa}-{\aa}}+(\tilde{\BB}- {\BB})\mm^2+\tilde{\BB}\hh^2.
	\]
	Hence,
	\begin{align*}    
		\hh&=\underbrace{\tilde{R}_{\mm}^{-1}(	{\tilde{\aa}-{\aa}}+(\tilde{\BB}- {\BB})\mm^2+\tilde{\BB}\hh^2)}_{:= \phi(\hh)}.
	\end{align*}
We wish to prove that, for a suitable choice of $k_{\hh}$, $\hh \in B_{k_{\hh}}$ implies $\phi(\hh) \in B_{k_{\hh}}$; i.e., $\phi(B_{k_{\hh}}) \subseteq B_{k_{\hh}}$. By the Brouwer fixed-point theorem, if this inclusion holds then there exists $\hh \in B_{k_{\hh}}$ such that $\hh = \phi(\hh)$, i.e., $\tilde{\mm} = \mm + \hh$ satisfies~\eqref{eq02}. 
 Since $\tilde{R}_{\mm}$ is an M-matrix, $\tilde{R}_{\mm}^{-1}\geq 0$.
Now, by applying  \cite[Theorem 2.5]{alfa}, we can derive that $
	\tilde{R}_{\mm}^{-1}\leq \gamma R_{\mm}^{-1}$. 
    Consequently, for $\hh \in B_{k_{\hh}}$, we obtain
	\begin{align*}
		\abs{\phi(\hh)}
		&\leq \tilde{R}_{\mm}^{-1}(|{\tilde{\aa}-{\aa}}|+|(\tilde{\BB}- {\BB})|\mm^2+\tilde{\BB}|\hh|^2)\\
        &\leq\gamma{R}_{\mm}^{-1}(\varepsilon\aa+\varepsilon{\BB}\mm^2+\tilde{\BB}|\hh|^2)\\
		&=\gamma{R}_{\mm}^{-1}(\varepsilon{\mm}+\tilde{\BB}|\hh|^2).
	\end{align*}
Form \eqref{kapax} and  $|\hh|\leq k_{\hh}\mm$, we have that
    \begin{align*}
		\abs{\phi(\hh)}
		\leq \gamma\left( \varepsilon(2\kappa-1)+(1+\varepsilon)(\kappa-1) k_{\hh}^2\right) \mm.
	\end{align*}
	If $k_{\hh}>0$ is chosen such that 
	\begin{equation} \label{condkh4}
		\gamma(2\kappa-1)\varepsilon+\gamma(\kappa-1)(1+\varepsilon) k_{\hh}^2= k_{\hh},    
	\end{equation}
	then the previous computation shows that $\phi(\hh) \in B_{k_{\hh}}$. It remains to show that a $k_{\hh}>0$ satisfying~\eqref{condkh4} exists.
	If we set $s = 1 / k_{\hh}$, \eqref{condkh4} becomes
	\begin{gather}
		\gamma(\kappa-1)(1+\varepsilon) k_{\hh}^2-k_{\hh}+\gamma(2\kappa-1)\varepsilon = 0, \nonumber \\
		\gamma(2\kappa-1)\varepsilon s^2-s+\gamma(\kappa-1)(1+\varepsilon)=0. \label{qp24}    
	\end{gather}
	The discriminant of this quadratic equation in $s$ is
	\[
	\Delta = 1-4\varepsilon\gamma^2(2\kappa-1)(\kappa-1)(1+\varepsilon),
	\]
in which $\Delta>0$ because of \eqref{condDelta},
hence the quadratic equation \eqref{qp24} has two real solutions
	$s_1,s_2=\frac{1 \pm \sqrt{\Delta}}{2\gamma(2\kappa-1)\varepsilon}$. Note that $2\gamma(2\kappa-1)\varepsilon> 0$ because of~$\kappa>1$. Hence, the largest solution $s_1$ satisfies
$s_{1}\geq\frac{1}{2\gamma(2\kappa-1)\varepsilon} > 0$.
It follows from $s={1}/{k_{\hh}}$ that $ k_{\hh}\leq 2\varepsilon(2\kappa-1)\gamma$, and thus $|\hh|\leq  (2\varepsilon(2\kappa-1)\gamma)\mm$, which completes the proof.
\end{proof}
This proof strategy with the Brouwer fixed-point theorem is classic; it appears for instance in~\cite[Theorem~V.2.11]{stewartsun}.
\subsection{Componentwise error bounds under zero-sum perturbations}
Let us now consider a perturbation of the multilinear PageRank problem \eqref{mlpr}, under the conditions $\mathbf{1}_n^T\tilde{\BB}_{(1)}=\mathbf{1}_n^T\BB_{(1)}=\alpha\mathbf{1}_{n^2}^T,$ and $\mathbf{1}_n^T\aa=\mathbf{1}_n^T\tilde{\aa}=1-\alpha$. Here all perturbations are zero-sum; In particular, for the perturbed vector $\tilde{\aa}$ and the perturbed tensor $\tilde{\BB} $, we can write $
\tilde{\aa}=\aa+ {\mathbf{e}}$ and $\tilde{\BB}={\BB}+{\mathscr{E}}
$,
where ${\mathscr{E}}$ and ${\mathbf{e}}$ satisfy $\mathbf{1}_n^T{\mathscr{E}}_{(1)}=\mathbf{0}$ and $\mathbf{1}_n^T{\mathbf{e}}=0$.

This kind of perturbation is natural in the context of applications of the multilinear PageRank, since $\vv$ and $\PP$ are taken normalized and hence guaranteed to be stochastic.

Consider a problem of the form~\eqref{mlpr} with $\alpha < 1/2$, and let $\mm$ be its minimal solution with $\mathbf{1}_n^T \mm = 1$. Then, $R_{\mm}^T$ is an M-matrix with $R_{\mm}^T \mathbf{1}_n = \mathbf{1}_n (1-2\alpha)$, and we can apply Theorem~\ref{thmRS} to it, constructing a decomposition $({R_{\mm}^T})^{-1} = R+S$.

For any vector $\hh$ such that $\mathbf{1}_n^T\hh = 0$, we have $S^T \hh = R_{\mm}^{-1}\hh$; hence, in our convergence theorems we replace $R_{\mm}^{-1}$ with $S^T$, which does not diverge when $\alpha \to 1/2$. To this purpose, we define $\mathbf{z}= S^T \mm$, and
\begin{equation}\label{omega}
\omega = \max_{m_i\neq 0} \frac{z_i}{m_i},    
\end{equation}
analogously to $\kappa$ in~\eqref{kapavalue}, so that $\mathbf{z} \leq \omega \mm$.

\begin{theorem}\label{thm_4}
	Let $\BB=(b_{ijk})$ and $\tilde{\BB}=(\tilde{b}_{ijk})$ be nonnegative $n$-dimensional tensors. Let $\mm$ be the minimal solution of~\eqref{eq01}. Suppose that $d(\tilde{\BB},\BB) \leq \varepsilon$
    and $d(\tilde{\aa},{\aa}) \leq \varepsilon$ 
    for a $\varepsilon$ small enough to ensure~\eqref{Mmatrixcond} and
	\begin{equation}\label{cond_Delta}
		\varepsilon+\varepsilon^2 <\frac{1}{4\gamma^2\omega^2},    
	\end{equation}  
	 where $\gamma=\frac{(1+\varepsilon)^{n-1}}{(1-\varepsilon)^{n}}$, $\omega$ is defined in \eqref{omega} and $S$ is the bounded matrix defined in Theorem \ref{thmRS}.
	We suppose that $\mathbf{1}_n^T\tilde{\BB}_{(1)}=\mathbf{1}_n^T\BB_{(1)}$ and $\mathbf{1}_n^T\aa=\mathbf{1}_n^T\tilde{\aa}$.
	Then, \eqref{eq11v} has a solution $\tilde{\mm}$ such that 
	\[
	d(\tilde{\mm},\mm)\leq 2\omega\gamma \varepsilon.
	\]
\end{theorem}
\begin{proof}
	We define the set $
	B_{k_{\hh}}=\left\{{\hh}\in \mathbb{R}^{n}: \mathbf{1}_n^T{\hh}=0,\,\abs{\hh} \leq k_{\hh} \mm\right\}$.
	For a vector $\hh \in B_{k_{\hh}}$, we have that $\tilde{\mm} = \mm + \hh$ satisfies~\eqref{eq11v} if and only if
	\begin{align*}
		\hh=\tilde{\mm}-\mm&={\tilde{\aa}-{\aa}}+\tilde{\BB}(\mm+\hh)^2-{\BB}\mm^2\\
		&={\tilde{\aa}-{\aa}}+ \tilde{\BB}\mm^2+(\tilde{\BB}\mm\hh+\tilde{\BB}\hh\mm)+\tilde{\BB}\hh^2-{\BB}\mm^2.	
	\end{align*}
	This relation can be rewritten as
	\[
	\tilde{R}_{\mm}\hh:={\hh-\tilde{\BB}\mm\hh-\tilde{\BB}\hh\mm} =
	{\tilde{\aa}-{\aa}}+(\tilde{\BB}- {\BB})\mm^2+\tilde{\BB}\hh^2.
	\]
	Hence,
	\begin{align*}    
		\hh&=\underbrace{\tilde{R}_{\mm}^{-1}(	{\tilde{\aa}-{\aa}}+(\tilde{\BB}- {\BB})\mm^2+\tilde{\BB}\hh^2)}_{:= \phi(\hh)}.
	\end{align*}
	We wish to prove that, for a suitable choice of $k_{\hh}$, $\hh \in B_{k_{\hh}}$ implies $\phi(\hh) \in B_{k_{\hh}}$; i.e., $\phi(B_{k_{\hh}}) \subseteq B_{k_{\hh}}$. By the Brouwer fixed-point theorem, this inclusion implies that there exists $\hh \in B_{k_{\hh}}$ such that $\hh = \phi(\hh)$, i.e., $\tilde{\mm} = \mm + \hh$ satisfies~\eqref{eq11v}. Since $\tilde{R}_{\mm}$ is an M-matrix, Theorem~\ref{thmRS} shows that $\tilde{S}^T\leq \gamma S^T$.
	Consequently, for $\hh \in B_{k_{\hh}}$, we obtain
	\begin{align*}
		|\phi(\hh)|
        &={|\tilde{S}^T({\tilde{\aa}-{\aa}}+(\tilde{\BB}- {\BB})\mm^2+\tilde{\BB}\hh^2)|}\\
		&\leq \tilde{S}^T(|{\tilde{\aa}-{\aa}}|+|(\tilde{\BB}- {\BB})|\mm^2+\tilde{\BB}|\hh|^2)\\
        &\leq\gamma{S}^T(\varepsilon\aa+\varepsilon{\BB}\mm^2+\tilde{\BB}|\hh|^2)\\
		&=\gamma{S}^T(\varepsilon{\mm}+\tilde{\BB}|\hh|^2).
	\end{align*}
	It now follows from $S^T\mm\leq \omega\mm$, $\tilde{\BB}\leq(1+\varepsilon){\BB}$ and $|\hh|\leq k_{\hh}\mm$ that
    \begin{align*}
		\abs{\phi(\hh)}\leq \gamma(\varepsilon\omega \mm +{S}^T\tilde{\BB}|\hh|^2)\leq \gamma(\varepsilon\omega \mm +(1+\varepsilon){S}^T{\BB}|\hh|^2)\leq \gamma \varepsilon\omega \mm+\gamma (1+\varepsilon)k_{\hh}^2S^T{\BB}\mm^2.
	\end{align*}
    Since $\mm=\aa+\BB\mm^2$ and $\mm,\aa,\BB\mm^2\geq \mathbf{0}$, we have $\BB\mm^2\leq \mm$. Hence,
    \begin{align*}
		\abs{\phi(\hh)}\leq \gamma \varepsilon\omega \mm+\gamma (1+\varepsilon)k_{\hh}^2S^T\mm\leq \gamma\varepsilon\omega\mm+\gamma(1+\varepsilon)k_{\hh}^2 \omega\mm.
	\end{align*}
	If $k_{\hh}>0$ is chosen such that, 
	\begin{equation} \label{cond_kh4}
		\gamma\omega(\varepsilon+(1+\varepsilon) k_{\hh}^2)= k_{\hh},    
	\end{equation}
	then the previous computation shows that $\phi(\hh) \in B_{k_{\hh}}$. It remains to show that such $k_{\hh}>0$ exists.
	If we set $s = 1 / k_{\hh}$, \eqref{cond_kh4} becomes
	\begin{gather}
		\gamma\omega(1+\varepsilon) k_{\hh}^2-k_{\hh}+\gamma\omega\varepsilon = 0, \nonumber \\
		\gamma\omega\varepsilon s^2-s+\gamma\omega(1+\varepsilon)=0. \label{qp_24}    
	\end{gather}
	The discriminant of this quadratic equation in $s$ is
	\[
	\Delta = 1-4\varepsilon\gamma^2\omega^2(1+\varepsilon),
	\]
	in which $\Delta>0$ because of \eqref{cond_Delta},
	hence the quadratic polynomial \eqref{qp_24} has two real solutions
	$s_1,s_2=\frac{1 \pm \sqrt{\Delta}}{2\gamma\omega\varepsilon}$. Note that $2\gamma\omega\varepsilon> 0$ because of~$\omega>0$. Hence, the largest solution $s_1$ satisfies
	$
	s_{1}\geq\frac{1}{2\gamma\omega\varepsilon} > 0.
	$
	It follows from $s={1}/{k_{\hh}}$ that $
	k_{\hh}\leq 2\omega\gamma\varepsilon$, 	and $|\hh|\leq  (2\varepsilon\omega\gamma)\mm$, which completes the proof.
\end{proof}
\section{GTH-based methods for the minimal solution}\label{sec:algorithms}
This section proposes some algorithms to achieve high componentwise accuracy for the solution of the multilinear PageRank problem \eqref{mlpr}.
\subsection{Newton-GTH}
The multilinear PageRank problem requires highly accurate solutions, especially in a componentwise sense. Although the standard Newton method is effective for many non-linear systems, it often struggles with the ill-conditioned Jacobian matrices that arise in this problem, leading to numerical instability and loss of precision in small components. To solve this problem, we propose the Newton-GTH method, a new approach that combines the Newton method with the GTH algorithm. Specifically, we use the GTH algorithm to solve the inner linear systems at each step of the Newton method, and use a subtraction-free formulation of the residual. This makes the Newton-GTH method effective in mitigating rounding errors and achieving componentwise accuracy. This method also maintains reliable performance even when the Jacobian matrices are close to being singular.

We recall that the GTH algorithm~\cite{alfa,Ocinneide} is a variant of Gaussian elimination that can be used to produce the LU factorization of an M-matrix $M$ in a subtraction-free fashion, hence computing a result with accuracy that is always of the order of machine precision, irrespective of the condition number $\kappa(M)$.

The GTH algorithm requires as inputs not the M-matrix $M$, but a so-called \emph{triplet representation}, i.e., a triple $(-\offdiag(M),\vv,\mathbf{w})$ such that $\vv>\mathbf{0}, \mathbf{w} \geq \mathbf{0}$ and $M\vv = \mathbf{w}$. (Note that the diagonal entries of $M$ can be reconstructed from the relation $M\vv = \mathbf{w}$.) Alternatively, one can apply the algorithm to the transpose of $M$, requiring then two w vectors such that $\vv^TM=\mathbf{w}^T$. This is the form that we shall use in the following.

We would like to write down the Newton method for equation \eqref{mlpr}. Given a first solution guess $\xx$, we would like to find $\hh$ so that $\xx+\hh$ is a better approximation of the solution. Using equation  \eqref{mlpr}, we obtain that
\begin{align*}
	{\xx+\hh} &=(1-\alpha){\vv}+\alpha\PP{({\xx+\hh})}^2\\
	&=(1-\alpha){\vv}+\alpha\PP\xx^2+\alpha\PP\xx\hh+\alpha\PP\hh\xx+\alpha\PP\hh^2. 
\end{align*}
If we omit the second-order term $\alpha\PP\hh^2$, this relation can be written as
\[
R_{\xx}\hh =(1-\alpha){\vv}+\alpha\PP\xx^2-\xx.
\]
Hence we can formulate the Newton method for this problem as 
\begin{equation}\label{eqNn}
 \begin{split}    
\hh & =  R_{\xx_k}^{-1}\left((1-\alpha){\vv}+\alpha\PP\xx_k^2-\xx_k\right), \\
\xx_{k+1} & =\xx_k + \hh.
 \end{split}
\end{equation}
In the case $\alpha < 1/2$, we can start the method from $\xx_0 = \mathbf{0}$. There are general results in~\cite{qve} that ensure that this method converges monotonically. Moreover, again in the case $\alpha < 1/2$, it is guaranteed that $R_{\xx_k}$ is a nonsingular M-matrix for each $k$. We would like to use the GTH algorithm to compute its inverse; for this, we need to compute a triplet representation for its transpose $R_{\xx_k}^T$ at each step $k$.
\begin{theorem}\label{thmz}
Consider the Newton method \eqref{eqNn}. Let $z_k = 1 - 2\alpha (\mathbf{1}_n^T \xx_k)$; then $\mathbf{1}_n^T R_{\xx_k} = z_k \mathbf{1}_n^T$ gives a triplet representation for $R_{\xx_k}$. Moreover, the scalars $z_k$ satisfy the recurrence relation
\begin{equation} \label{zrecurrence}
z_{k+1} = \frac{(1 - 2\alpha)^2 + z_{k}^2}{2z_{k}}, \quad k \geq 0,
\end{equation}
with the initial point $z_0 = 1 - 2\alpha (\mathbf{1}_n^T \xx_0) = 1$. In particular, if the value of $1-2\alpha$ is known with high accuracy, we can compute each $z_k$ without cancellation.
\end{theorem}
\begin{proof}
For the first part, we have to verify that
\begin{align*}
\mathbf{1}_n^TR_{\xx_k}=\mathbf{1}_n^T(I_n-\alpha\PP{{:}\xx_k}-\alpha\PP{\xx_k{:}})
=(1-2\alpha (\mathbf{1}_n^T\xx_k) )\mathbf{1}_n^T.
\end{align*}
Now we show that the scalars $z_k$ can be computed with a recursion, starting from $z_0 = 1- 2\alpha(\mathbf{1}_n^T\xx_0) = 1$. 
Multiplying both sides of~\eqref{mlpr} by $\mathbf{1}_n^T$,  we obtain
\begin{equation*}
(1-2\alpha (\mathbf{1}_n^T\xx_k))\mathbf{1}_n^T\hh= 1-\alpha+\alpha(\mathbf{1}_n^T\xx_k)^2-\mathbf{1}_n^T\xx_k. 
\end{equation*}
Hence,
\begin{equation*}
\mathbf{1}_n^T\hh=\frac{1-\alpha+\alpha(\mathbf{1}_n^T\xx_k)^2-\mathbf{1}_n^T\xx_k}{1-2\alpha (\mathbf{1}_n^T\xx_k)}.
\end{equation*}
We can get
\begin{equation*}
{1}-2\alpha (\mathbf{1}_n^T\xx_{k+1})=1-2\alpha(\mathbf{1}_n^T\xx_k)-2\alpha(\mathbf{1}_n^T\hh)=1-2\alpha(\mathbf{1}_n^T\xx_k)-2\alpha(\frac{1-\alpha+\alpha(\mathbf{1}_n^T\xx_k)^2-\mathbf{1}_n^T\xx_k}{1-2\alpha (\mathbf{1}_n^T\xx_k)}).
\end{equation*}	
This simplifies to
\begin{equation*}
	1-2\alpha (\mathbf{1}_n^T\xx_{k+1})=
	\frac{1}{2}\frac{(1-2\alpha)^2+(1-2\alpha (\mathbf{1}_n^T\xx_k))^2}{1-2\alpha (\mathbf{1}_n^T\xx_k)},
\end{equation*}	
or equivalently, $z_{k+1}=\frac{(1-2\alpha)^2+z_k^2}{2z_k}$. This completes the proof.
\end{proof}
Another point we need to obtain a subtraction-free implementation of the Newton method is computing after each step the residual vector $\mathbf{r}_{k+1}$; this vector can be written in a subtraction-free form as
\begin{align*}
\mathbf{r}_{k+1}&=(1-\alpha)\vv+\alpha\PP{\xx^2_{k+1}}-\xx_{k+1}\\
&=(1-\alpha)\vv+\alpha\PP{(\xx_{k}+\hh)}^2-(\xx_{k}+\hh)\\
&=(1-\alpha)\vv+\alpha\PP\xx_{k}^2+\alpha\PP\xx_{k}\hh+\alpha\PP\hh\xx_{k}+\alpha\PP\hh^2-\xx_{k}-\hh\\
&=(1-\alpha)\vv+\alpha\PP\xx_{k}^2-\xx_{k}-R_{\xx_k}\hh+\alpha\PP\hh^2\\
&=\alpha\PP\hh^2.
\end{align*}
Algorithm \ref{algo1} summarizes the subtraction-free implementation of the Newton-GTH method as discussed above.
	\begin{algorithm}[h]
			\caption{Newton method based on GTH (Newton-GTH).}
			\label{algo1}
			\begin{algorithmic}
					\STATE{\textbf{Input}: A stochastic vector $\vv$, a stochastic tensor $\PP$, a parameter $\alpha \in (0,1)$, a tolerance $\varepsilon>0$.}
					\STATE{\textbf{Output}: An approximation to the minimal solution $\mm$ of \eqref{mlpr}.}
					\STATE{$\xx  \gets \mathbf{0}$; $z \gets 1$};
					\COMMENT{Invariant: at each step, $z = 1-2\alpha(\mathbf{1}_n^T\xx)$}
		           \STATE{$\mathbf{r}\leftarrow (1-\alpha)\vv$;}
		            \COMMENT{Invariant: at each step, $\mathbf{r}=(1-\alpha)\vv+\alpha\PP\xx^2-\xx$}
					\WHILE{${\|\mathbf{r}\|}_{\infty} > \varepsilon$}
                \STATE{Compute $L,U$ factors such that $LU=R_{\xx}$, using the GTH algorithm on the triplet $(\offdiag(\alpha\PP{{:}\xx}+\alpha\PP{\xx{:}}),\mathbf{1}_n^T,z\mathbf{1}_n^T)$;} 
					\STATE{  $\hh\leftarrow U^{-1} (L^{-1}\mathbf{r})$;}
					\STATE{$\xx\leftarrow\xx+\hh$;} 	    \STATE{$z\leftarrow\frac{(1-2\alpha)^2+z^2}{2z}$;}
					\STATE{$\mathbf{r}\leftarrow \alpha\PP\hh^2$;}
					\ENDWHILE
                      \STATE{$\mm  = \xx$;}
				\end{algorithmic}
		\end{algorithm}
\subsection{The block Jacobi-GTH}\label{secBJ}
In this section, we propose a block Jacobi-type method for solving the multilinear PageRank problem \eqref{mlpr}. The need to improve on the Newton method arises from the observation that in large dimension solving a linear system with $R_{\xx_k}$ might be challenging, since even when $\PP$ is sparse the matrix $R_{\xx_k}^{-1}$ might be dense, as well as the LU factors or $R_{\xx_k}$. To avoid this problem, we define a new method which introduces an M-matrix splitting $R_{\xx_k} = M_{\xx_k} - N_{\xx_k}$, and replaces the inversion with one single iteration of a classical iterative method, starting from the initial guess $\xx_k$.

The most simple such iterative method is the block Jacobi method; for more details on this iteration, the reader is referred to \cite[Subsection 4.1.1]{saad} and \cite[Section 8.5]{stoer}. Using the second version of the update in~\eqref{newton}, the Newton iteration for the PageRank problem \eqref{mlpr} can be written as
\begin{align}\label{erx}
R_{\xx_{k}}\xx_{k+1} =(1-\alpha){\vv}-\alpha\PP\xx_{k}^2,
\end{align}
where $R_{\xx_{k}}=I_n-\alpha\PP{{:}\xx_{k}}-\alpha\PP{\xx_{k}{:}}$. We assume that $n=ps$, and write~\eqref{erx} as a block-structured linear system
\begin{equation}\label{eqk1}
	R_{\xx_{k}}\xx_{k+1}\equiv
	\begin{bmatrix}
		R_{11}&\dots& R_{1s}\\
		\vdots&\ddots&\vdots\\
		R_{s1}&\dots& R_{ss}
	\end{bmatrix}
	\begin{bmatrix}
		\xx_{k+1}^{(1)}\\\vdots \\\xx_{k+1}^{(s)}
	\end{bmatrix}=\begin{bmatrix}
		((1-\alpha){\vv}-\alpha\PP\xx_{k}^2)^{(1)}\\\vdots \\((1-\alpha){\vv}-\alpha\PP\xx_{k}^2)^{(s)}
	\end{bmatrix},
\end{equation}
where $R_{\xx}=[R_{ij}]_{i,j=1}^s$ is a block-structured matrix with $s$ square diagonal blocks $R_{ii}$ of size $ p \times p$. 
We consider the splitting $ R_{\xx_{k}}=M_{\xx_{k}}-N_{\xx_{k}}$, where 
\[M_{\xx_{k}}=
\begin{bmatrix}
	R_{11}&&&O\\
	&	R_{22}&&\\
	&&\ddots&\\
	O&&& R_{ss}
\end{bmatrix},\quad N_{\xx_{k}}=-
\begin{bmatrix}
	O&R_{12}&\dots&R_{1s}\\
	R_{21}&\ddots&\ddots&\vdots\\
	\vdots&\ddots&\ddots&R_{s-1,s}\\
	R_{s1}&\dots&R_{s,s-1}&O
\end{bmatrix}.
\]
One step of the block Jacobi iteration for~\eqref{eqk1} then corresponds to
\begin{equation*}
	R_{ii}{({\xx}_{k+1})}^{(i)}={\left({N_{\xx_k}\xx_{k}}\right)}^{(i)}+{\left({(1-\alpha){\vv}-\alpha\PP\xx_{k}^2}\right)}^{(i)}, \quad i=1,2,\dots,s,\quad k=0,1,\dots,
\end{equation*}
or, equivalently,
\begin{equation}\label{bljac}
M_{\xx_{k}}\xx_{k+1} = N_{\xx_k}\xx_k + (1-\alpha){\vv}-\alpha\PP\xx_{k}^2.
\end{equation}
We can rewrite the above block Jacobi iteration as follows
 \begin{align}
 \xx_{k+1} &= M_{\xx_k}^{-1}((M_{\xx_k} - R_{\xx_k})\xx_k +  (1-\alpha){\vv}-\alpha\PP\xx_{k}^2) \nonumber\\
 &= \xx_k + M_{\xx_k}^{-1}( (1-\alpha){\vv}+\alpha\PP\xx_{k}^2-\xx_k) \nonumber\\
 &= \xx_k + M_{\xx_k}^{-1}F(\xx_k).
 \label{eqJN}
 \end{align}
This form highlights the similarity of the block Jacobi method with the Newton method: both are fixed-point methods with a similar form, with the difference that at each iteration of the block Jacobi we solve a linear system not with the matrix $R_{\xx_k}$, but with $M_{\xx_k}$, which is block diagonal; hence the solution can be computed more efficiently.

\begin{remark}
    In a similar way, one could consider applying an iterative method to the first formulation of the Newton method~\eqref{newton}: to solve the linear system $\hh = \xx_{k+1} - \xx_k = R_{\xx_k}^{-1}F(\xx_k)$, we replace the inversion with one step of an iterative method based on an M-matrix splitting. However, unlike in~\eqref{erx} where we could use $\xx_{k}$ as a starting point, there is no natural starting point to use for $\hh$ in this formulation.
\end{remark}
\begin{remark}
    In~\cite[Section~4]{qve}, another family of methods based on a splitting is considered; they are different from the ones studied here, since in~\cite{qve} we had (after switching to the notation in this paper) $M-N = I_n - \BB\xx_k^2 \neq R_{\xx_k}$.
\end{remark}

\subsubsection{Triplet representation}
To implement the GTH algorithm within the block Jacobi method \eqref{bljac}, we establish the following triplet representation.

Let $R_{\yy_k} = M_{\yy_k} - N_{\yy_k}$ be the splitting for the block Jacobi method, where $R_{\yy_k}=I_n-\alpha \PP\yy_k{:}- \alpha \PP{:}\yy_k$.
	Consider the sequence $\yy_k$ generated by the block Jacobi iteration:
	\begin{equation}\label{BJy}
		M_{\yy_k}\yy_{k+1} = N_{\yy_k}\yy_k + (1-\alpha){\vv}-\alpha \PP\yy_{k}^2.
	\end{equation}
In the case $\alpha < 1/2$, we can start the method from $\yy_0 = \mathbf{0}$.  
Define $u_k = 1- 2\alpha (\mathbf{1}_n^T\yy_k)$. This gives $\mathbf{1}_n^T R_{\yy_k} = u_k \mathbf{1}_n^T$. Since $R_{\yy_k} = M_{\yy_k} - N_{\yy_k}$, then
\[\mathbf{1}_n^T M_{\yy_k}=\mathbf{1}_n^T N_{\yy_k}+u_k\mathbf{1}_n^T,\]
is a triplet representation for the block Jacobi method.
 We shall show that the scalars $u_k$ can be computed with a recursion, starting from $u_0 = 1- 2\alpha(\mathbf{1}_n^T\yy_0) = 1$. 
	 Subtracting $N_{\yy_k}{\yy}_{k+1}$ from both sides of \eqref{BJy} yields 
	\begin{equation*}
		M_{\yy_k}{\yy}_{k+1}-N_{\yy_k}{\yy}_{k+1}+N_{\yy_k}({\yy}_{k+1}-{\yy}_{k})=(1-\alpha)\vv-\alpha \PP\yy_k^2.
	\end{equation*}
	Multiplying by $\mathbf{1}_n^T$ from the left gives
	\begin{equation*}
		(1-2\alpha (\mathbf{1}_n^T\yy_k))\mathbf{1}_n^T{\yy}_{k+1}= 1-\alpha-\alpha(\mathbf{1}_n^T\yy_k)^2-\mathbf{1}_n^TN_{\yy_k}({\yy}_{k+1}-{\yy}_{k}).
	\end{equation*}
	This leads to
	\begin{equation*}
		\mathbf{1}_n^T{\yy}_{k+1}=\frac{1-\alpha-\alpha(\mathbf{1}_n^T\yy_k)^2-\mathbf{1}_n^TN_{\yy_k}({\yy}_{k+1}-{\yy}_{k})}{1-2\alpha (\mathbf{1}_n^T\yy_k)}.
	\end{equation*}
	Consequently, we can derive
	\begin{equation*}
		{1}-2\alpha (\mathbf{1}_n^T\yy_{k+1})=1-2\alpha(\frac{1-\alpha-\alpha(\mathbf{1}_n^T\yy_k)^2-\mathbf{1}_n^TN_{\yy_k}({\yy}_{k+1}-{\yy}_{k})}{1-2\alpha (\mathbf{1}_n^T\yy_k)}).
	\end{equation*}	
	This simplifies to
	\begin{equation*}
		1-2\alpha (\mathbf{1}_n^T\yy_{k+1})=
		\frac{(1- 2\alpha(\mathbf{1}_n^T\yy_k))^2+(1-2\alpha)^2+4\alpha\mathbf{1}_n^TN_{\yy_k}({\yy}_{k+1}-{\yy}_{k})}{2(1-2\alpha (\mathbf{1}_n^T\yy_k))}.
	\end{equation*}	
	Equivalently,
	\begin{equation} \label{urecurrence}
	u_{k+1}=\frac{u_k^2+(1-2\alpha)^2+4\alpha \mathbf{1}_n^TN_{\yy_k}({\yy}_{k+1}-{\yy}_{k})}{2u_k},\quad k\geq 0.	    
	\end{equation} 
	Here, $u_{k}$ defines the triplet representation for the block Jacobi method, i.e.,
	\begin{equation}\label{jtriplet}
		R_{\yy_k} = \left(\offdiag(\alpha(\PP{:}\yy_k + \PP{\yy_k{:}})), \mathbf{1}_n^T, u_{k}\mathbf{1}_n^T\right).
	\end{equation}
 Or, equivalently,
\begin{equation*}\label{jtriplet1}
M_{\yy_k} = \left(-\offdiag(M_{\yy_k}), \mathbf{1}_n^T,\mathbf{1}_n^TN_{\yy_k}+u_{k}\mathbf{1}_n^T\right).
\end{equation*} 
The following theorem establishes the monotonic convergence of the block Jacobi method. 
\begin{theorem}\label{thmj}
Consider the splitting $R_{\xx_k} = M_{\xx_k} - N_{\xx_k}$, where $R_{\xx_k} = I_n - \BB\xx_k{:} - \BB{:}\xx_k$ is an M-matrix. The block Jacobi iteration is defined as
	\begin{equation}\label{eqbj}  
		M_{\xx_k} \xx_{k+1} = N_{\xx_k} \xx_k + \aa- \BB \xx_k^2.  
	\end{equation}
	Assume that $\mathbf{0}\leq \xx_k\leq \xx_*$ and $F(\xx_k)=\aa+\BB \xx_k^2-\xx_k\geq \mathbf{0}$ for $k\geq 0$. 
	Then, $\xx_k\leq \xx_{k+1}\leq \xx_*$, and $F(\xx_{k+1})=\BB\hh^2+N_{\xx_{k}}\hh\geq \mathbf{0}$, where $\hh = \xx_{k+1} - \xx_{k}$. In particular, the generated sequence  $\{\xx_k\}$ converges monotonically to $\xx_{*}$.
\end{theorem}
\begin{proof}
From \eqref{eqJN}, we have $\xx_{k+1}=\xx_{k}+M_{\xx_{k}}^{-1}F(\xx_k)$. Since $M_{\xx_{k}}^{-1}\geq 0$ and $F(\xx_k)\geq \mathbf{0}$, it follows that $\xx_{k+1}-\xx_{k}\geq \mathbf{0}$. Now, we show that $F(\xx_{k+1})\geq \mathbf{0}$. We rewrite \eqref{eqbj} as
	\[  
	R_{\xx_{k}} \xx_{k+1} = (M_{\xx_{k}} - N_{\xx_{k}}) \xx_{k+1} = N_{\xx_{k}} (\xx_{k} - \xx_{k+1}) + \aa- \BB \xx_{k}^2.  
	\]  
	This simplifies to
	\[  
	\xx_{k+1} - \BB(\xx_{k} \xx_{k+1} + \xx_{k+1} \xx_{k}) = \aa- \BB \xx_{k}^2 - N_{\xx_{k}} (\xx_{k} - \xx_{k+1}).  
	\]  
	Hence,  
	\begin{align*}  
		F(\xx_{k+1}) &= \aa+ \BB \xx_{k+1}^2 - \xx_{k+1} \\ 
		&= \aa+ \BB \xx_{k+1}^2 - \xx_{k+1}+\left(\xx_{k+1} - \BB(\xx_{k} \xx_{k+1} + \xx_{k+1} \xx_{k}) - \aa+ \BB \xx_{k}^2 + N_{\xx_{k}} (\xx_{k} - \xx_{k+1})\right)\\
		&= \BB\hh^2 + N_{\xx_{k}} \hh \geq \mathbf{0}.
	\end{align*}  
	 Subsequently, we show that $\xx_{k+1} \leq \xx_* $.
	\begin{align*}  
		R_{\xx_{k}} (\xx_* - \xx_{k+1}) &= (M_{\xx_{k}} - N_{\xx_{k}})(\xx_* - \xx_{k+1}) \\  
		&= \xx_* - \BB(\xx_{k} \xx_* + \xx_* \xx_{k}) - N_{\xx_{k}} (\xx_{k} - \xx_{k+1}) - \aa+ \BB \xx_{k}^2 \\  
		&= N_{\xx_{k}} \hh + \BB(\xx_* - \xx_{k})^2\geq \mathbf{0}.  
	\end{align*}  
	Since $R_{\xx_{k}}^{-1}\geq 0$, it follows that $\xx_* \geq \xx_{k+1} $. The sequence $\xx_k$ is monotonic and bounded from above by $\xx_*$, so it converges.
\end{proof}

The convergence of the Newton and block Jacobi methods is compared in the following theorem.
\begin{theorem}\label{thmcompareBJN}
	Let $F(\xx) = \aa + \BB \xx^2 - \xx$, and suppose that $R_{\xx} = I_n - \BB\xx{:} - \BB{:}\xx$, where $\xx \in [\mathbf{0},\mm]$ and $\mm$ is the minimal solution. Let $R_{\yy_k} = M_{\yy_k} - N_{\yy_k}$ be the splitting for the block Jacobi method. Consider the sequences $\{\xx_k\}$ (Newton method) and $\{\yy_k\}$ (block Jacobi method) generated as follows:
	\begin{align*}
		\xx_{k+1} &= \xx_k + R_{\xx_k}^{-1}F(\xx_k), \\
		\yy_{k+1} &= M_{\yy_k}^{-1}(N_{\yy_k}\yy_k + \aa - \BB  \yy_k^2).
	\end{align*}
 Assume $\mathbf{0}\leq \xx_k \leq \mm$, $\mathbf{0}\leq \yy_k \leq \mm$ and $\yy_k\leq \xx_k$. Suppose that  $F(\xx_k)\geq \mathbf{0}$ and $F(\yy_k) \geq \mathbf{0}$. Then $\xx_{k+1}\geq \yy_{k+1}$,  for all $k \geq 0$. Moreover, $F(\xx_{k+1})\geq \mathbf{0}$ and $F(\yy_{k+1}) \geq \mathbf{0}$.
\end{theorem}
\begin{proof}
From \eqref{eqJN}, the block Jacobi iteration can be written as $\yy_{k+1}=\yy_k+M_{\yy_k}^{-1}F(\yy_k)$.
Define $\zz=\xx_k+R_{\yy_k}^{-1}F(\xx_k)$.
We now prove that $\zz-\yy_{k+1}\geq \mathbf{0}$. We get
\begin{align*}
\zz-\yy_{k+1}=\xx_{k}+R_{\yy_{k}}^{-1}F(\xx_k)-\yy_k-M_{\yy_k}^{-1}F(\yy_k).
\end{align*}	
It follows from $R_{\yy_k}\leq M_{\yy_k}$ that
$R_{\yy_k}^{-1}\geq M_{\yy_k}^{-1} \geq 0$. Consequently, we obtain,
\begin{align*}
	\zz-\yy_{k+1}&\geq \xx_k-\yy_k+R_{\yy_{k}}^{-1}(F(\xx_k)-F(\yy_k))\\
	&=R_{\yy_{k}}^{-1}(R_{\yy_{k}}\xx_k-R_{\yy_{k}}\yy_k+F(\xx_k)-F(\yy_k))\\
    &=R_{\yy_{k}}^{-1}(\xx_{k}-\BB\xx_{k}\yy_{k}-\BB\yy_{k}\xx_{k}-\yy_{k}+2\BB\yy_{k}^2+\aa+\BB\xx_{k}^2-\xx_{k}-\aa-\BB\yy_{k}^2+\yy_{k})\\
    &=R_{\yy_{k}}^{-1}(-\BB\xx_{k}\yy_{k}-\BB\yy_{k}\xx_{k}+\BB\yy_{k}^2+\BB\xx_{k}^2)\\
	&=R_{\yy_{k}}^{-1}\BB(\xx_{k}-\yy_{k})^2\geq \mathbf{0}.
\end{align*}	
We now show that $\xx_{k+1}-\zz\geq \mathbf{0}$.
Since $\xx_k \geq \yy_k$, we have $
\BB\xx_k{:} + \BB{:}\xx_k \geq \BB\yy_k{:} + \BB{:}\yy_k$,
which implies $R_{\xx_k} \leq R_{\yy_k}$, and therefore $0\leq R_{\yy_k}^{-1} \leq R_{\xx_k}^{-1}$. 
Since $F(\xx_k) \geq \mathbf{0}$, we obtain
\[
\xx_{k+1} - \zz = \xx_k + R_{\xx_k}^{-1} F(\xx_k) - \xx_k - R_{\yy_k}^{-1} F(\xx_k)= \bigl( R_{\xx_k}^{-1} - R_{\yy_k}^{-1} \bigr) F(\xx_k) \geq \mathbf{0}.
\]
From $\yy_{k+1}\leq \zz$ and $\zz\leq \xx_{k+1}$, we conclude that $\yy_{k+1}\leq \xx_{k+1}$. By Theorem \ref{thm.21}, we have $F(\xx_{k+1})\geq \mathbf{0}$. In addition, Theorem \ref{thmj} ensures that $F(\yy_{k+1})\geq \mathbf{0}$. Thus, the proof is complete.
\end{proof}
This theorem shows that the block Jacobi method always converges slower than the Newton method, when started from the same vector, for instance $\xx_0 = \yy_0 = \mathbf{0}$. 
In particular, one also has 
\[u_k = 1- 2\alpha (\mathbf{1}_n^T\yy_k)\geq 1- 2\alpha (\mathbf{1}_n^T\xx_k)=z_k.\]
This observation motivates a variant of the block Jacobi method: we can use at each step instead of $u_k$ the formula for $z_k$ from the Newton method, dropping the summand $4\alpha \mathbf{1}_n^TN_{\yy_k}({\yy}_{k+1}-{\yy}_{k})$ in~\eqref{urecurrence}.

    \subsubsection*{Variant of the block Jacobi-GTH method}
    We modify~\eqref{bljac} and define a new numerical method to solve~\eqref{mlpr} as
    \[
    T_k \ww_{k+1}= N_{\ww_k}\ww_k + (1-\alpha)\vv - \alpha \PP \ww_k^2,
    \]
    where at each step $T_k$ is the M-matrix with triplet representation
	\begin{equation*}
		T_k = \left(\offdiag(-M_{\ww_k}), \mathbf{1}_n^T, \mathbf{1}_n^T N_{\ww_k} + z_k\mathbf{1}_n^T\right),
	\end{equation*}
	where $z_{k+1}=\frac{z_k^2 + (1-2\alpha)^2}{2z_k}$. 
    Hence, the only difference between the block Jacobi-GTH and its variant lies in the triplet representation: we use~\eqref{zrecurrence} from the Newton method, rather than~\eqref{urecurrence} from block Jacobi. The full algorithm is presented in Algorithm \ref{algo2}.
\begin{algorithm}[h]
		\caption{The block Jacobi variant based on GTH algorithm to solve \eqref{mlpr}.}
		\label{algo2}
		\begin{algorithmic}
			\STATE \textbf{Input}:
             A stochastic vector $\vv$, a stochastic tensor $\PP$, a parameter $\alpha \in (0,1)$, a tolerance $\varepsilon>0$, a dimension $p$ of each block diagonal, where $sp=n$ and $s$ is a positive integer.
			\STATE{\textbf{Output}: An approximation to a minimal solution $\mm$ of \eqref{mlpr}.}
			\STATE{$\ww  \gets \mathbf{0}$; $z \gets 1$};
			\WHILE{${\|\mathbf{r}\|}_{\infty} > \varepsilon$}
            \STATE{$\mathbf{b} \leftarrow (1-\alpha)\vv - \alpha \PP \ww^2$;}
            \STATE{$R\leftarrow   I_n-\alpha\PP{{:}\ww}-\alpha\PP{\ww{:}}$; 
            Compute $N$ from  $R$;}
            \STATE $\mathbf{z}\leftarrow \mathbf{1}_n^TN+z\mathbf{1}_n^T
			$; 
            \STATE{$\ww_{\text{old}} \leftarrow \ww$;}
			\FOR{$i=1,\dots,s$}
              \STATE{Compute $L,U$ factors such that $LU=(T_k)_{ii}$ using the GTH algorithm on the triplet $(-\offdiag(R_{ii}),\mathbf{1}^T_p,\mathbf{z}^{(i)})$};
			\STATE  $\ww^{(i)}\leftarrow U^{-1}\,\left(L^{-1}\, ((N_{(i),:}\ww_{\text{old}})^{(i)} +\mathbf{b}^{(i)})\right)$;
			\ENDFOR
            \STATE{$z\leftarrow\frac{(1-2\alpha)^2+z^2}{2z}$;}
		 	\ENDWHILE
                \STATE{$\mm  = \ww$;}
		\end{algorithmic}
	\end{algorithm}	
    Note that $T_k$ is guaranteed to be a nonsingular M-matrix at each step, since it has a triplet representation with $z_k \neq 0$.

    This algorithm has lower computational cost than the Newton method, since at each step we need to solve $s$ small linear system of sizes summing up to $n$, rather than a large one of size $n$.

    At least in the initial iterations, the scalar $z_k$ is strictly smaller than $1-2\alpha(\mathbf{1}^T\ww_k)$; hence $T_k \leq M_{\ww_k}$ and $T_k^{-1} \geq M_{\ww_k}^{-1}$. This means that this variant takes longer steps towards the minimal solution $\mm$. Unfortunately, these longer steps also mean that we cannot guarantee $\ww_k \leq \mm$: an analogue of~Theorem~\ref{thmcompareBJN} cannot be proved for this variant, and indeed experimentally one can observe examples in which $\ww_k \leq \mm$ fails to hold. Hence monotonic convergence properties do not hold for this variant. However, as we shall see in the experiments, on problems with $\alpha$ close to $1/2$ this variant displays in practice a convergence speed that is much closer to that of the Newton method than to that of the block-Jacobi method. So we can obtain the cheaper cost per iteration of block Jacobi together with the faster convergence speed of Newton, with only a small change to the triplet.

\section{Limiting accuracy of the Newton method for $\alpha > 1/2$} \label{sec:limitaccuracy}

For the case $\alpha > 1/2$, unfortunately it seems that very little can be said about componentwise accuracy: in the Newton method, $R_{\xx_*}^{-1}$ has mixed signs, and when working with a stochastic vector $\mathbb{x}$ mixed signs seem impossible to avoid in the computation of $F(\xx) = \xx-\aa-\BB\xx^2$, since $\mathbf{1}_n^T F(\xx) = 0$.

However, we can use the sign structure of our equation to improve our bounds on the \emph{limiting accuracy} of the Newton method in floating-point arithmetic, i.e., the value at which $\norm{\xx_{k}-\xx_*}$ stagnates. The limiting accuracy of the multivariate Newton method is studied in~\cite{tisseur}; the results obtained there (Theorem~2.2, Corollary~2.3) show that one can expect stagnation when $\norm{\xx_{k+1}-\xx_*} \approx \norm{R_{\xx_*}^{-1}}\norm{\mathbf{e}}$, where $\norm{\mathbf{e}}$ bounds the error incurred when computing the residual $\mathbf{r}_k$ (in our case $\mathbf{r}_k = \xx_k - \aa - \BB \xx_k^2$) in floating-point arithmetic.

We give here a slightly modified version of those results, to show that $\norm{R_{\xx_*}^{-1}}\norm{\mathbf{e}}$ can be replaced by $\norm{R_{\xx_*}^{-1}\mathbf{e}}$.
\begin{theorem}[\protect{Variant of~\cite[Theorem~2.2]{tisseur}}]
Let $\tilde{\xx}_k$ be the sequence computed by Newton's method for~\eqref{eq01} in machine arithmetic with precision $\mathsf{u}$. 

Define:
\begin{itemize}
    \item $\mathbf{e}_k$ the error incurred in the computation of the residual $\tilde{\mathbf{r}}_k$;
    \item $E_k$ the error incurred when computing the Jacobian $R_{\tilde{\xx}_k}$ and solving the linear system $R_{\tilde{\xx}_k}^{-1}(\mathbf{r}_k+\mathbf{e}_k)$, counted as a backward error for the linear system, so that the increment computed in machine arithmetic is 
    \[
    \tilde{\hh}_k = (R_{\tilde{\xx}_k}+E_k)^{-1}(\mathbf{r}_k+\mathbf{e}_k);
    \]
    \item $\boldsymbol{\varepsilon}_k$ the error incurred in the subtraction $\tilde{\xx}_k - \tilde{\hh}_k$.
\end{itemize}
Suppose that
\[
\norm{R_{\tilde{\xx}_k}^{-1}E_k} \leq \nu < 1, \quad \beta \norm{R_{\xx_*}^{-1}} \norm{\tilde{\xx}_k - \xx_*} \leq \mu < 1,
\]
where $\beta$ is a Lipschitz constant for $R_{\xx}$. (We can take $\beta= \frac{\norm{\BB}}{2}$ in our problem.) Then,
\[
\norm{\tilde{\xx}_{k+1} - \xx_*} \leq G_k \norm{\tilde{\xx}_{k} - \xx_*} + \mathbf{g}_k,
\]
with
\[
G_k = \frac{\nu}{1-\nu} + \frac{\mu}{2(1-\mu)(1-\nu)}  + \mathsf{u}, \quad \mathbf{g}_k = \frac{1}{1-\nu}\frac{1}{1-\mu}\norm{R_{\xx_*}^{-1}\mathbf{e}_k} + \mathsf{u}(\norm{\xx_*} + \norm{\tilde{\hh}_k}).
\]
\end{theorem}
\begin{proof}
    We have
    \[
        \tilde{\xx}_{k+1} - \xx_* = \tilde{\xx}_k - \xx_* - (R_{\tilde{\xx}_k}+E_k)^{-1}(\mathbf{r}_k+\mathbf{e}_k) + \boldsymbol{\varepsilon}_k,
    \]
    which gives
    \begin{multline*}
    \norm{\tilde{\xx}_{k+1} - \xx_*} \leq \norm{I_n-(R_{\tilde{\xx}_k}+E_k)^{-1}R_{\tilde{\xx}_k}}\norm{\tilde{\xx}_k - \xx_*} + \norm{(R_{\tilde{\xx}_k}+E_k)^{-1}}\norm{\mathbf{r}_k - R_{\tilde{\xx}_k}(\tilde{\xx}_k - \xx_*)} \\ {}+ \norm{(R_{\tilde{\xx}_k}+E_k)^{-1}\mathbf{e}_k } + \norm{\boldsymbol{\varepsilon}_k}
    \end{multline*}
    We bound three of these four summands using the same arguments as in~\cite[Proof of Theorem~2.2]{tisseur}:
    \begin{align*}
        \norm{I_n-(R_{\tilde{\xx}_k}+E_k)^{-1}R_{\tilde{\xx}_k}}\norm{\tilde{\xx}_k - \xx_*} &\leq \frac{\nu}{1-\nu} \norm{\tilde{\xx}_k - \xx_*}\\
        \norm{(R_{\tilde{\xx}_k}+E_k)^{-1}}\norm{\mathbf{r}_k - R_{\tilde{\xx}_k}(\tilde{\xx}_k - \xx_*)} & \leq \frac{1}{(1-\mu)(1-\nu)} \norm{R_{\xx_*}^{-1}} \cdot \frac{\beta}{2} \norm{\tilde{\xx}_k - \xx_*}^2\\
        \norm{\boldsymbol{\varepsilon}_k} &\leq \mathsf{u}(\norm{\tilde{\xx}_k - \xx_*} + \norm{\xx_*} + \norm{\tilde{\hh}_k}).
    \end{align*}
     For the remaining term, we need to modify that argument, to get
    \begin{align*}
        \norm{(R_{\tilde{\xx}_k}+E_k)^{-1}\mathbf{e}_k } &=
        \norm{(I_n+R_{\tilde{\xx}_k}^{-1}E_k)^{-1}R_{\tilde{\xx}_k}^{-1}\mathbf{e}_k }\\
        &\leq \norm{(I_n+R_{\tilde{\xx}_k}^{-1}E_k)^{-1}}\norm{R_{\tilde{\xx}_k}^{-1}\mathbf{e}_k}\\
        &\leq \frac{1}{1-\nu}\norm{R_{\tilde{\xx}_k}^{-1}\mathbf{e}_k}\\
        &\leq \frac{1}{1-\nu}\norm{(I_n+R_{\xx_*}^{-1}(R_{\tilde{\xx}_k}-R_{\xx_*}))^{-1}}\norm{R_{\xx_*}^{-1}\mathbf{e}_k}\\
        &\leq \frac{1}{1-\nu}\frac{1}{1-\mu}\norm{R_{\xx_*}^{-1}\mathbf{e}_k}.
    \end{align*}
\end{proof}
In particular, if $G_k$ and $\mathbf{g}_k$ are bounded uniformly by $G\leq \frac12$ and $\mathbf{g}$, respectively, from a certain $k$ onwards, we expect the iteration to converge in machine arithmetic, and stagnate at a limit accuracy of $\mathbf{g}$. Hence the key factor for the limiting accuracy is the error in the computation of the residual $\mathbf{e}_k$. We give a result to bound this error; to state it, we use the standard notation in machine error analysis $\tilde{\gamma}_m = \frac{Cm\mathsf{u}}{1-Cm\mathsf{u}}$ for a moderate constant $C$~\cite{higham-accuracy}.
\begin{theorem}
    When $\tilde{\xx}_k - \xx_* = o(1)$, we have
    \begin{equation} \label{ekbound}
    \abs{\mathbf{e_k}} \leq \tilde{\gamma}_{n^2} \xx_* + o(\mathsf{u}),  
    \end{equation}
    and hence
    \[
    \norm{R_{\xx_*}^{-1}\mathbf{e}_k}\leq \tilde{\gamma}_{n^2} \norm{\abs{R_{\xx_*}^{-1}}\xx_*} + o(\mathsf{u}).
    \]
\end{theorem}
\begin{proof}
Let $\mathbf{p}_k = \aa + \BB \tilde{\xx}_k^2$, and $\tilde{\mathbf{p}}_k$ its approximation computed in machine arithmetic starting from $\aa, \BB$ and~$\tilde{\xx}_k$. Since there is no cancellation in these sums of nonnegative terms, we have by standard arguments in error analysis
\[
\abs{\tilde{\mathbf{p}}_k - \mathbf{p}_k} \leq \tilde{\gamma}_{n^2} \mathbf{p}_k,
\]
 (Note, in passing, that the factor $n^2$ can be replaced with the maximum number of nonzero entries in a row of $\BB$.)

The error in the final subtraction $\tilde{\xx}_k - \tilde{\mathbf{p}}_k$ can be bounded by
\[
\abs{\tilde{\mathbf{r}}_k - (\tilde{\xx}_k - \tilde{\mathbf{p}}_k)} \leq \mathsf{u} \abs{\tilde{\xx}_k - \tilde{\mathbf{p}}_k}.
\]
Putting together these two bounds, we get
\[
\abs{\mathbf{e}_k} = \abs{\tilde{\mathbf{r}}_k - (\tilde{\xx}_k - \mathbf{p}_k)} \leq \tilde{\gamma}_{n^2} \mathbf{p}_k + \mathsf{u}  \abs{\tilde{\xx}_k - \tilde{\mathbf{p}}_k}.
\]
When $\tilde{\xx}_k - \xx_* = o(1)$ we have also $\tilde{\xx}_k - \mathbf{p}_k = o(1)$ and
$\tilde{\xx}_k - \tilde{\mathbf{p}}_k = (\tilde{\xx}_k - {\mathbf{p}}_k) - (\tilde{\mathbf{p}}_k - \mathbf{p}_k) = o(1)$, which gives~\eqref{ekbound}.
\end{proof}
This bound ultimately comes from the fact that the $\mathbf{r}_k$ is the subtraction of two positive quantities which can be computed with componentwise accuracy; it shows that the limiting accuracy of Newton's method is of the order of $\mathsf{u}\norm{\abs{R_{\xx_*}^{-1}}\xx_*}$, which might be significantly smaller than $\mathsf{u}\norm{R_{\xx_*}^{-1}}\norm{\xx_*}$, as shown in our experiments.

\section{Experimental results}\label{sec:experiments}\
In this section, we present some numerical experiments to assess the performance of componentwise accurate methods as well as componentwise perturbation bounds. All numerical experiments were computed using \textsc{Matlab} R2024a on a computer equipped with a 64-bit Intel Core i7-13620H processor. 

We use these abbreviations: the Newton method (\textbf{N}), the Newton-GTH method (\textbf{NG}, Algorithm \ref{algo1}), the block Jacobi-GTH variant (\textbf{BJGv}, Algorithm \ref{algo2}) and the block Jacobi method (\textbf{BJ}, Algorithm \ref{algo2}, but using \textsc{Matlab}’s backslash instead of GTH).

In the figures, $\ec$ and $\en$ represent the componentwise and normwise errors:
\[\ec=d(\tilde{\xx},\xx)=\max\limits_{i=1,\dots,n}\frac{|{x}_i-{x}^*_i|}{{x}^*_i},\qquad\en=\frac{{\|\xx-\xx^*\|}_2}{{\|\xx^*\|}_2},\]
where $\xx^*$ is the exact solution (either the minimal or a stochastic one) computed using variable-precision arithmetic (vpa).

For Example \ref{ex3}, we report results after 100 iterations, whereas in the other examples, results are provided after 500 iterations. All values of $\alpha$ are chosen so that an exact solution can be computed, achieving a residual error of approximately $10^{-40}$. In all examples, for the stochastic solution with $\alpha>1/2$,  the initial guess $\xx_0$ is set to $\vv$; otherwise, it is set to $\mathbf{0}$.

In the block Jacobi-GTH method, the dimension of the inner systems must be divisible by the block size $p$. In Examples \ref{ex1} and \ref{ex2}, we set $p=2$; for Example \ref{ex3}, the value of $p$ varies depending on the dataset and is given in its description.
\begin{example}\cite[Section 6.5]{gleich2015} \label{ex3}\rm
We examine an application in real-world networks. We consider an undirected graph $G$ containing $n$ nodes and construct a tensor $\mathscr{C}\in\mathbb{R}^{n\times n \times n}$ where
 \[ \mathscr{C}_{ijk}=\begin{cases}
1,&\text{there is a directed three-cycle between nodes}\,\, i,j,k,\\
0,&\text{otherwise}.
 \end{cases}	\]
We transform $\mathscr{C}$ into a stochastic tensor $\PP$. To do this, we first obtain ${S}\in\mathbb{R}^{n\times n^2}$ which is constructed from $\mathscr{C}$ through column-wise normalization
	\begin{equation*}
		S_{:,i} = 
		\begin{cases}
{(\mathscr{C}_{(1)})}_{:,i}/\mathbf{1}_n^T{(\mathscr{C}_{(1)})}_{:,i} & \text{if } \mathbf{1}_n^T{(\mathscr{C}_{(1)})}_{:,i} \neq 0, \\
			{(\mathscr{C}_{(1)})}_{:,i} & \text{otherwise}.
		\end{cases}
	\end{equation*}
 In particular, $S$ is a sub-stochastic matrix.
Let $A$ represent the adjacency matrix of the \texttt{ibm32} dataset \cite{Davis}. Let $D$ be the diagonal matrix with
 $d_{ii}=\sum_{j}A_{ij}$. We define $M = A^TD^\dagger$, where  $D^\dagger$ denotes the Moore-Penrose pseudo-inverse of $D$. For a matrix $B\in \mathbb{R}^{n\times q}$, we let $\mathrm{dangling}(B):=\mathbf{1}_q^T-\mathbf{1}_n^TB$. The tensor $\PP$ in \eqref{mlpr} is given by
 \[
 \PP_{(1)}=\nu({S}+\vv\,\,\mathrm{dangling}({S}))+(1-\nu)(M+\vv\,\,\mathrm{dangling}(M)) \otimes \mathbf{1}_n^T,
 \]
 where $\nu=0.1$ and $\vv = \mathrm{rand}(n{,}1).\exp(9  \mathrm{randn}(n{,}1))$. We normalize $\vv$ so that $\mathbf{1}_n^T\vv=1$; the vector $\vv$ is constructed in this way to ensure that the computed solution exhibits varying magnitudes across its components. The dimension $p$ of each block diagonal in \textbf{BJGv} for \texttt{ibm32} is set to 2.

Figure \ref{fig1} shows the numerical results for the minimal solution. 
Figure \ref{fig1a} demonstrates that the  \textbf{BJGv} and \textbf{NG} methods achieve high componentwise accuracy, whereas the \textbf{N} method stagnates at around $10^{-7}$. This is due to the ill-conditioning of $R_{\xx_k}$ when $\xx_k \approx \mm$, which prevents one from computing the Newton increment $\hh$ accurately. The use of the GTH algorithm circumvents this problem, since it allows to compute $R_{\xx_k}^{-1}$ in a componentwise accurate fashion.

Further comparison results for the computation of the minimal solution using the \textbf{BJGv} and \textbf{NG} methods are presented in Table \ref{T1} when other datasets are applied. The results are reported when ${\|\mathbf{r}_k\|}_{\infty}\leq 10^{-15}$ or a maximum of 100 iterations are achieved. In \textbf{BJGv}, the dimensions of each diagonal block $p$ is
12 for \texttt{Erdos991},
10 for \texttt{Harvard500},
7 for \texttt{Edinburgh}, and
14 for \texttt{NotreDame\_yeast}, respectively.
The table also shows that the  \textbf{BJGv}
requires more iterations than the \textbf{NG} to satisfy the stopping criterion ${\|\mathbf{r}_k\|}_{\infty}\leq 10^{-15}$ . However,
it achieves this with less CPU time compared to the \textbf{NG}.
\begin{figure}[H]
	\centering
	\centerline{\input{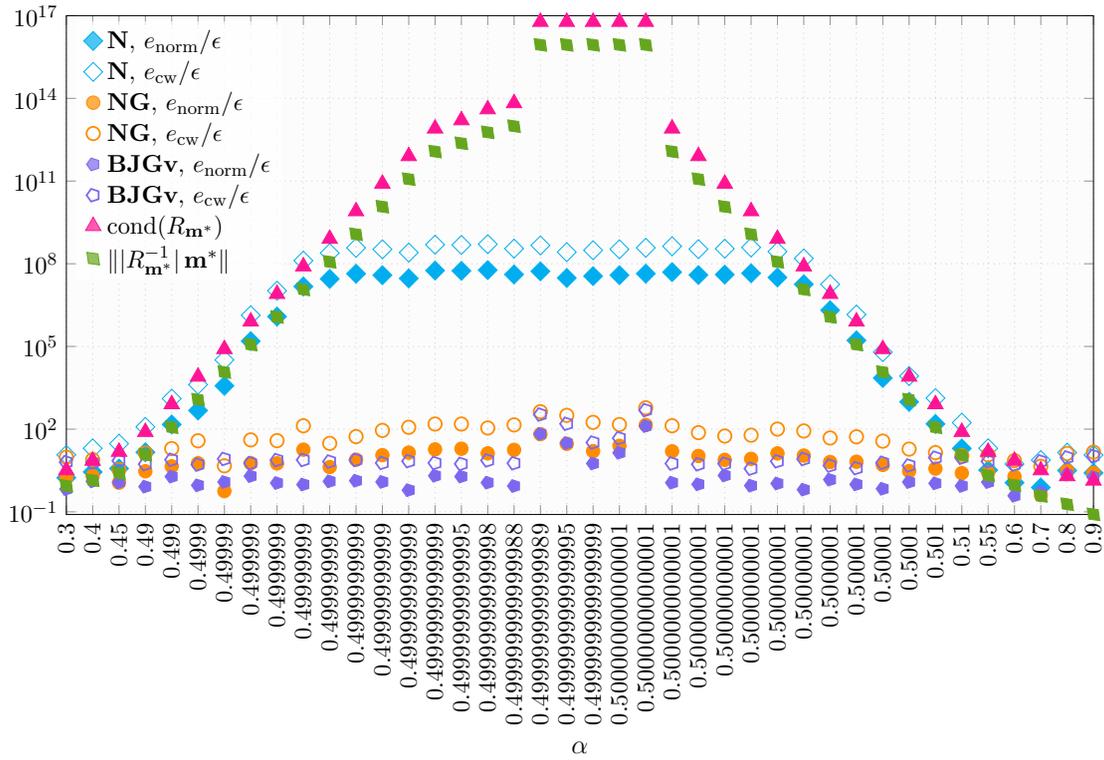}}
	\caption{Example \ref{ex3}:  Results obtained from the computation of the minimal solution $\mm$.
	}
    	\label{fig1}
\end{figure}
\begin{figure}[H]
	\centering
	\centerline{\input{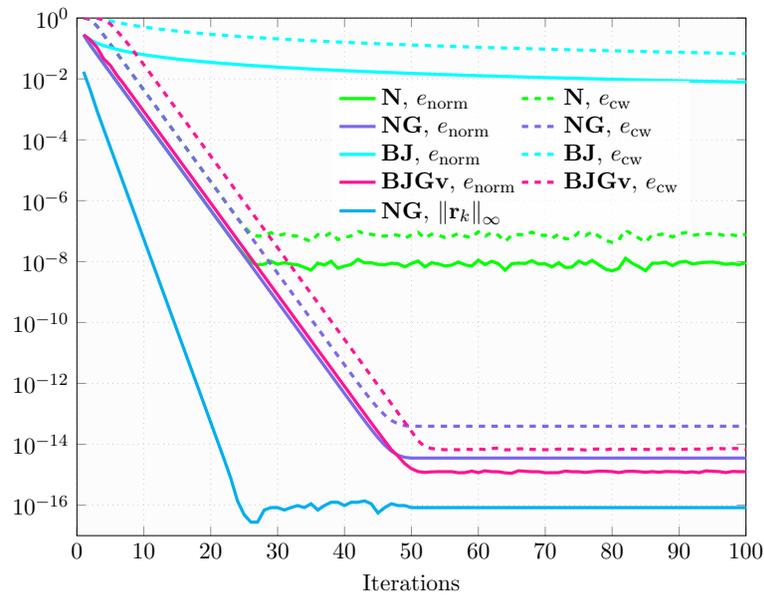}}
	\caption{Example \ref{ex3}:  Results for the minimal solution $\mm$ with $\alpha=0.499999999999999$. The residual for the \textbf{NG} method is also shown here (not invariant).}
    	\label{fig1a}
\end{figure}
\begin{table}[htbp]
	\centering
	\caption{Example \ref{ex3}: Performance comparison between \textbf{BJGv} and \textbf{NG} methods for different datasets to obtain the minimal solution $\mm$ where $\alpha=0.499999999999999$. ``IT'' denotes the number of iterations, and ``CPU'' denotes the CPU time in seconds.}
	\label{T1}
	\begin{tabular}{l l S[table-format=2.0] S[table-format=3.5] S[table-format=1.4e-2]}
		\toprule
		Dataset & Method & {IT} & {CPU} & ${\|\mathbf{r}_k\|}_{\infty}$ \\
		\midrule
        \multirow{2}{*}{\texttt{Erdos991} \cite{batage}} 
		& \textbf{BJGv} & 49 & {0.56} &9.0536e-16\\
		& \textbf{NG}  & 24 & {2.72} & 8.4152e-16 \\
        \cmidrule(r){1-5}
        \multirow{2}{*}{\texttt{Harvard500} \cite{moler}} 
		& \textbf{BJGv} & 48 & {0.58} & 7.4344e-16 \\
		& \textbf{NG}  & 24 & {2.91} & 4.1754e-16\\
        \cmidrule(r){1-5}
		\multirow{2}{*}{\texttt{Edinburgh} \cite{grindrod}} 
		& \textbf{BJGv} & 50 & {5.18} & 6.7704e-16 \\
		& \textbf{NG}  & 25 & {138.76} & 4.4514e-16\\
		\cmidrule(r){1-5}
		\multirow{2}{*}{\texttt{NotreDame\_yeast} \cite{Davis}} 
		& \textbf{BJGv} & 51 & {8.11} & 6.8665e-16\\
		& \textbf{NG}  & 24 & {297.64} &8.1545e-16 \\
		\bottomrule
	\end{tabular}
\end{table}
\end{example}
\begin{example}\label{ex1}\rm
We examine PageRank \eqref{mlpr} where the tensor $\PP$ and the vector $\vv$ are given by
\[
	\PP_{(1)}= 
	\left[\begin{array}{cccc:cccc:cccc:cccc}
		0&0&0&0&0&0&1&0&0&1&1&0&0&0.5&0&1\\
		0&0&0&0&0&0&0&0&0&0&0&0.5&1&0&0&0\\
		0&0&0&0&0&0&0&0&0&0&0&0.5&0&0&0&0\\
		1&1&1&1&1&1&0&1&1&0&0&0&0&0.5&1&0
	\end{array}\right],\,
    \vv=
	\begin{bmatrix}
	  1.5462\cdot 10^{-2}\\
	1.4317\cdot 10^{-12}\\
	3.5898\cdot 10^{-7}\\
	9.8454\cdot 10^{-1}
	\end{bmatrix}.
\]
For all selected $\alpha$, $\mm^\ast$ has entries with different orders of magnitude; for instance, for $\alpha=0.49999$, the exact solution is
$ \mm^\ast=
\begin{bmatrix}
    2.4655\cdot 10^{-1},
   8.2687\cdot 10^{-2},
   2.1565\cdot 10^{-7},
   6.7076\cdot 10^{-1}
   \end{bmatrix}^T$.
Figure \ref{fi} shows the numerical results for the minimal solution. 
\begin{figure}[h]
	\centering
	\centerline{\input{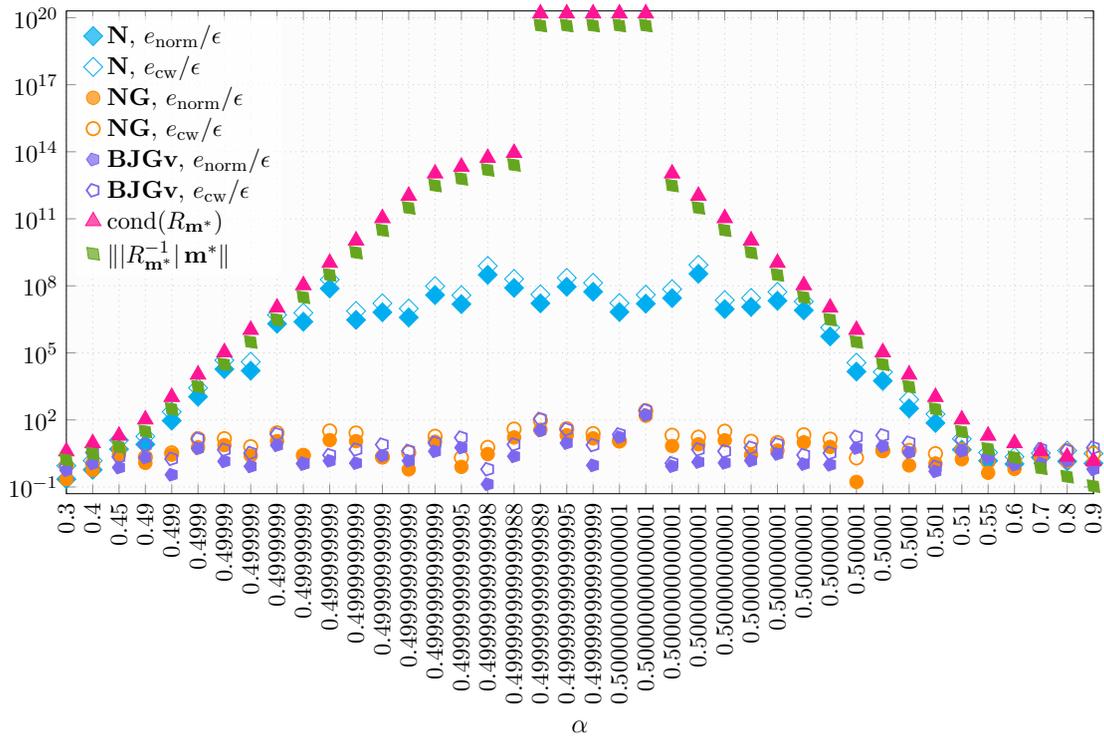}}
	\caption{Example \ref{ex1}:  Results obtained from the computation of the minimal solution $\mm$.
	}
	\label{fi}
\end{figure}
\end{example}
\begin{example}\label{ex2}\rm
We examine PageRank \eqref{mlpr} where the tensor $\PP$ is given by
\[
		\PP_{(1)}= 
		\left[\begin{array}{cccc:cccc:cccc:cccc}
0&0&0&0&0&0&0.5&0&1&0&0&0&0&0&0&0\\
0&0&0&0&0&0&0&1&0&1&0&0&0.5&0&0&0\\
0&0&0&0&0&0&0.5&0&0&0&1&0&0.5&0.5&1&0\\
1&1&1&1&1&1&0&0&0&0&0&1&0&0.5&0&1
		\end{array}\right],
		\]
  and  $\vv=  \begin{bmatrix}
			 10^{-4},
			0,
			0,
			9.999\cdot 10^{-4}
		\end{bmatrix}^T$.
 For all selected $\alpha$, the exact solution has entries with different orders of magnitude; e.g., for $\alpha=0.9951$, the exact solution is
 \[\ss^\ast=
	\begin{bmatrix}
      8.6225\cdot 10^{-7},
   8.5301\cdot 10^{-5},
   8.5971\cdot 10^{-3},
   9.9132\cdot 10^{-1}
   \end{bmatrix}^T.\]
Figure \ref{fig_ex2} shows that an increase in the condition number of the Jacobian matrix leads to larger componentwise and normwise errors.
\begin{figure}[H]
	\centering
	\centerline{\definecolor{marker1}{RGB}{27,158,119}  
\definecolor{marker3}{RGB}{117,112,179} 
\definecolor{marker4}{RGB}{231,41,138} 
\definecolor{marker2}{RGB}{255,140,0}       
\definecolor{marker5}{RGB}{124,104,238}    

 \begin{tikzpicture}[scale=.88]
 	\begin{axis}[%
 		width=5.5in,
 		height=3in,
 		scale only axis,
 		xmin=1,
 		xmax=23,
 		xtick={1,2,3,4,5,6,7,8,9,10,11,12,13,14,15,16,17,18,19,20,21,22,23},
 		xticklabels={
 			0.7,
 			0.8,
 			0.9,
 			0.92,
 			0.95,
 			0.97,
 			0.99,
 			0.993,
 			0.995,
 			0.9951,
 			0.99515,
 			0.99517,
 			0.995174,
 			0.9951747,
 			0.99517473,
 			0.995174734,
 			0.9951747342,
 			0.99517473429,
 			0.995174734291,
 			0.9951747342915,
 			0.99517473429157,
 			0.995174734291573,
 			0.9951747342915737
 		},
 		xticklabel style={rotate=90},
 		xlabel={$\alpha$},
 		xlabel style={shift={(0,2.5cm)},font=\large},
 		ymode=log,
 		ymin=0.009,
 		ymax=1e12,
 		legend style={
 			at={(0.22,0.996)}, 
 			anchor=north east, 
 			legend cell align=left, 
 			align=left,
 			draw=none,
 		},
 		grid=major,
 		grid style={dotted,gray!50},
 		axis background/.style={fill=gray!2},
 		only marks,
 		]
 		\addplot [color=marker5,
 		mark=square*,
 		mark size=2.5pt]
 		table[row sep=crcr]{%
 			1	0.562851686202694\\
 			2	0.0151361639920011\\
 			3	0.298861024497538\\
 			4	0.602917031736688\\
 			5	0.0905862123850469\\
 			6	0.820026083269263\\
 			7	0.290772774837082\\
 			8	0.209453740590997\\
 			9	1.3073869648\\
 			10	9.93450141729207\\
 			11	120.832962361125\\
 			12	160.639351823464\\
 			13	267.015611005961\\
 			14	261.412409812024\\
 			15	17481.9187201155\\
 			16	114641.320980319\\
 			17	327987.575126572\\
 			18	2056250.63434815\\
 			19	301775.689149764\\
 			20	99815.0706461064\\
 			21	4219992.36754801\\
 			22	63359628.9645749\\
 			23	20411151.1207587\\
 		};
 		\addlegendentry{\textbf{N}, $\en/ \epsilon$}
 		
 		\addplot [	color=marker2,
 		mark=square,
 		mark size=2.5pt,
 		mark options={solid, draw=marker5, fill=white, line width=0.8pt}]
 		table[row sep=crcr]{%
 			1	5.93381115047934\\
 			2	4.81022491758114\\
 			3	5.62372356601544\\
 			4	23.3246286542846\\
 			5	35.3146933619066\\
 			6	31.1158214583368\\
 			7	167.786210605385\\
 			8	27.5271790467104\\
 			9	160.478006279952\\
 			10	3857.71011300382\\
 			11	9288.06946108916\\
 			12	12517.3656674639\\
 			13	19221.4827021838\\
 			14	18076.2801336002\\
 			15	1249607.09703739\\
 			16	8183811.59618463\\
 			17	23419707.5278054\\
 			18	146830655.107386\\
 			19	21548459.622137\\
 			20	7128984.54885756\\
 			21	301335315.108475\\
 			22	4524283906.65971\\
 			23	1457487565.69879\\
 		};
 		\addlegendentry{\textbf{N}, $\ec/ \epsilon$}
 		
 		\addplot [	color=cyan,
 		mark=*,
 		mark size=2.5pt,
 		]
 		table[row sep=crcr]{%
 			1	0.562851647923736\\
 			2	0.0151588422953661\\
 			3	0.298878467953268\\
 			4	0.602909621778593\\
 			5	0.0901745047624065\\
 			6	0.819678916067999\\
 			7	0.0972095980888455\\
 			8	0.458240844725699\\
 			9	7.67683870852382\\
 			10	8.70929946914935\\
 			11	2.66762421895778\\
 			12	26.1552820752686\\
 			13	32.1043446283464\\
 			14	122.716220299655\\
 			15	1199.35945318889\\
 			16	574.644272263807\\
 			17	4014.93873969392\\
 			18	55221.3745739401\\
 			19	46235.8606210277\\
 			20	107495.570806715\\
 			21	403559.327103481\\
 			22	1089489.59649634\\
 			23	2386082.43440801\\
 		};
 		\addlegendentry{\textbf{NG}, $\en/ \epsilon$}
 		
 		\addplot [	color=cyan,
 		mark=o,
 		mark size=2.5pt,
 		mark options={solid, draw=cyan, fill=white, line width=0.8pt}
 		]
 		table[row sep=crcr]{%
 			1	2.06933400455826\\
 			2	8.37213308902808\\
 			3	12.9972005601547\\
 			4	16.7594431078319\\
 			5	17.5683861766649\\
 			6	12.4086418605682\\
 			7	98.7393600361617\\
 			8	78.2313238444404\\
 			9	752.788304258579\\
 			10	794.251540880521\\
 			11	460.405189910352\\
 			12	1837.08660455084\\
 			13	2255.20928658141\\
 			14	8743.80654537382\\
 			15	85627.5907774282\\
 			16	40958.2485856375\\
 			17	286473.54927806\\
 			18	3943184.45995724\\
 			19	3301492.22163718\\
 			20	7675870.41998507\\
 			21	28816646.7979705\\
 			22	77796445.4618231\\
 			23	170381391.450513\\
 		};
 		\addlegendentry{\textbf{NG}, $\ec/ \epsilon$}
 		
 		\addplot [
 			color=marker3,
 		mark=triangle*,
 		mark size=4.5pt,
 		mark options={rotate=90}
 		]
 		table[row sep=crcr]{%
 			1	21.1822814125324\\
 			2	39.357940385038\\
 			3	149.978489909296\\
 			4	235.017874993437\\
 			5	609.210325390836\\
 			6	1722.74174904729\\
 			7	17556.6323704716\\
 			8	43169.8540409298\\
 			9	231782.657649016\\
 			10	362878.257975245\\
 			11	637978.741244722\\
 			12	1464062.70302746\\
 			13	3719227.46482911\\
 			14	17208485.8479298\\
 			15	48641645.0680732\\
 			16	186609365.115053\\
 			17	332982094.888236\\
 			18	2540101437.78779\\
 			19	4207016234.44944\\
 			20	11743468157.8234\\
 			21	52531267286.1221\\
 			22	122035910915.388\\
 			23	283138837871.902\\
 		};
\addlegendentry{$\mathrm{cond}(R_{\ss^\ast})$}
 		
 		\addplot [
 			color=marker4,
 		mark=diamond*,
 		mark size=4.5pt,
 		mark options={rotate=45}
 		]
 		table[row sep=crcr]{%
 			1	2.49938320349423\\
 			2	1.66604941214037\\
 			3	1.24910873471165\\
 			4	1.18947097214075\\
 			5	1.11018851390814\\
 			6	1.06715448904302\\
 			7	1.5403136767879\\
 			8	3.58731680057707\\
 			9	26.6136927342626\\
 			10	44.1205643885347\\
 			11	81.570347515869\\
 			12	195.06904297948\\
 			13	505.734523724843\\
 			14	2364.92411998329\\
 			15	6697.4649874226\\
 			16	25714.0948750745\\
 			17	45889.2258582518\\
 			18	350105.208746649\\
 			19	579862.66457548\\
 			20	1618641.62541042\\
 			21	7240585.21722558\\
 			22	16820685.0070077\\
 			23	39026138.2843785\\
 		};
 		\addlegendentry{$\||R_{\ss^\ast}^{-1}|\,\ss^\ast\|$}   
        
      \addplot [
 			color=pink,
 		mark=diamond*,
 		mark size=4.5pt,
 		mark options={rotate=90}
 		]
  table[row sep=crcr]{%
1	15.043807136904\\
2	24.944367806585\\
3	84.3722179379276\\
4	129.110093067041\\
5	323.0501630065\\
6	892.385429436183\\
7	8880.33748509312\\
8	21747.3538300638\\
9	116310.0542436\\
10	182014.080234537\\
11	319889.561058926\\
12	733907.166829313\\
13	1864153.60637323\\
14	8624719.60825961\\
15	24378436.1914102\\
16	93525299.5066441\\
17	166884597.879294\\
18	1273051891.99108\\
19	2108478683.05894\\
20	5885608659.34325\\
21	26327697366.4114\\
22	61162136119.8585\\
23	141903936221.755\\
};
\addlegendentry{$\|R_{\ss^\ast}^{-1}\|\|\ss^\ast\|$}   

 	\end{axis}
 	
 \end{tikzpicture}%
}
	\caption{Example \ref{ex2}: Results for the computation of the stochastic solution $\ss$.
	}
	\label{fig_ex2}
\end{figure}
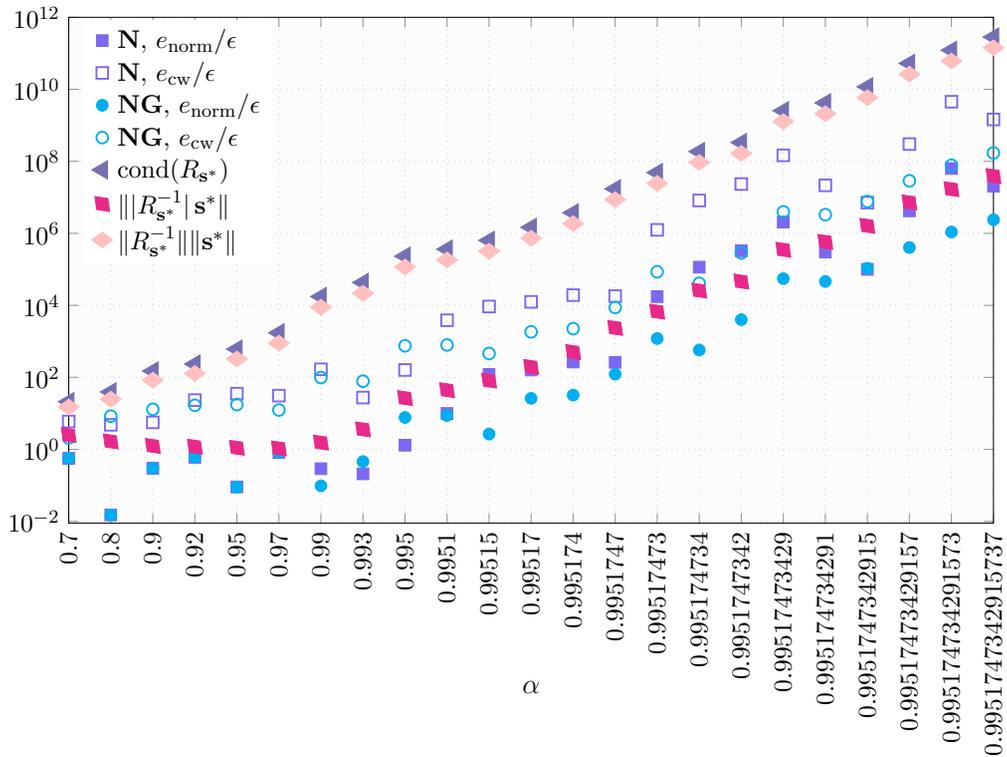
\end{example}
It is interesting to note that $\||R_{\ss^\ast}^{-1}|\,\ss^\ast\|$ is a much better predictor of the actual growth in the error than $\mathrm{cond}(R_{\ss^\ast})$ or $\|R_{\ss^\ast}^{-1}\|\|\ss^\ast\|$; this result confirms the value of our bound on the limiting accuracy.

\begin{example}\label{ex4}\rm
We consider equation \eqref{mlpr} where the tensor ${\mathscr{P}}$ and the vector $\vv$ are taken from Examples \ref{ex1} and \ref{ex2}. Let
$	\tilde{\vv}=\mathbf{v}+\varepsilon {\mathbf{e}}$ and $\tilde{\mathscr{P}}={\mathscr{P}}+\varepsilon{\mathscr{E}}$, where the perturbations ${\mathscr{E}}$ and ${\mathbf{e}}$ are generated by 
${\mathscr{E}}_{(1)} = {R} - \frac{1}{n}\mathbf{1}_n.\mathbf{1}_n^T{R},$ and $
{\mathbf{e}} =\mathbf{r} - \frac{1}{n}\mathbf{1}_n.\mathbf{1}_n^T\mathbf{r}$.
Here, ${R}=\mathrm{rand}(n{,}n^2)$ and $\mathbf{r}=\mathrm{rand}(n{,}1)$. These perturbations satisfy $\mathbf{1}_n^T{\mathscr{E}}_{(1)}=\mathbf{0}$  and $\mathbf{1}_n^T{\mathbf{e}}=0$. 
We set $\tilde{\PP}=\max(\tilde{\PP},0)$ and $
\tilde{\vv}=\max(\tilde{\vv},0)$ to ensure non-negativity.
Finally, we renormalize the perturbed tensor and vector with $ \tilde{\PP}_{(1)}=\tilde{\PP}_{(1)}./{\mathbf{1}_n^T\tilde{\PP}_{(1)}}, $ and $\tilde{\vv}=\tilde{\vv}/{\mathbf{1}_n^T\tilde{\vv}}$.
For all selected $\varepsilon$, we have ${\|\mathbf{1}_n^T\tilde{\PP}_{(1)}-\mathbf{1}_n^T{\PP}_{(1)}\|}_{\infty}\approx 10^{-16}$ and ${\|\mathbf{1}_n^T\tilde{\vv}-\mathbf{1}_n^T\vv\|}_{\infty}\approx 10^{-16}$. 

In Figures \ref{figure6}, \ref{figure7} and \ref{figure8}, we present the componentwise errors and their corresponding upper bounds derived from Theorem \ref{thm_4} for both minimal and stochastic solutions. Here, the stochastic solution with $\alpha > 1/2$ was computed using the Newton method, while the minimal solution was obtained by Newton-GTH.
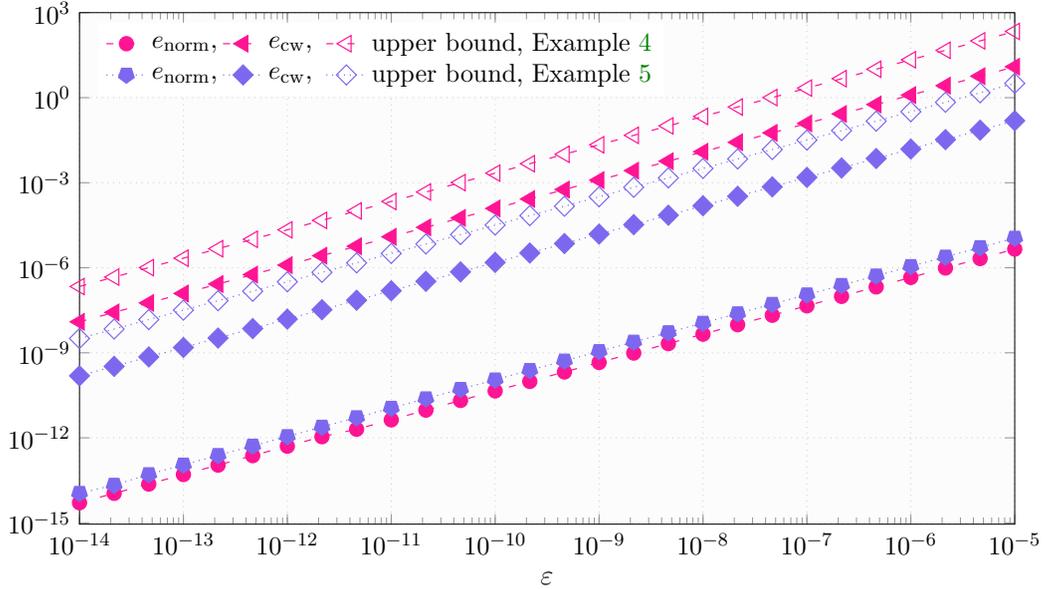
\begin{figure}[h]
 	\centerline{    
\definecolor{mycolor2}{RGB}{124,104,238}    
\definecolor{mycolor1}{RGB}{255,20,147}    
  \begin{tikzpicture}[scale=.89]
	
\begin{axis}[%
		width=5.5in,
		height=3in,
scale only axis,
xmode=log,
xmin=1e-14,
xmax=0.00001,
xminorticks=true,
ymode=log,
ymin=9e-16,
ymax=1000,
yminorticks=true,
		xlabel={$\varepsilon$},
				xlabel style={font=\large}, 
		legend style={
			at={(0.02,0.98)},
			anchor=north west,
			legend cell align=left,
			align=left,
			draw=none,
			font=\normalsize,
		},
            grid=major,
 		grid style={dotted,gray!50},
 		axis background/.style={fill=gray!2},
		legend columns=3
		]
\addplot [color=mycolor1, loosely dashed, line width=.5pt, mark=*, mark options={solid, fill=mycolor1, mycolor1,rotate=90},mark size=3pt]
  table[row sep=crcr]{%
1e-14	5.12597565875631e-15\\
2.1544e-14	1.10074041612996e-14\\
4.6416e-14	2.3424178517669e-14\\
1e-13	5.07954986707473e-14\\
2.1544e-13	1.09179218441504e-13\\
4.6416e-13	2.35117644963804e-13\\
1e-12	5.07839660560911e-13\\
2.1544e-12	1.09471596185079e-12\\
4.6416e-12	2.01536569766644e-12\\
1e-11	4.31126062502766e-12\\
2.1544e-11	9.50366035735578e-12\\
4.6416e-11	2.07951707281484e-11\\
1e-10	4.51726628547873e-11\\
2.1544e-10	9.77166439982151e-11\\
4.6416e-10	2.10923728338999e-10\\
1e-09	4.54828012081793e-10\\
2.1544e-09	9.80288985024993e-10\\
4.6416e-09	2.11241720511948e-09\\
1e-08	4.55146152983831e-09\\
2.1544e-08	9.80607280848711e-09\\
4.6416e-08	2.11273473856623e-08\\
1e-07	4.5517787903777e-08\\
2.1544e-07	9.8063920701214e-08\\
4.6416e-07	2.11276595044482e-07\\
1e-06	4.55180770921434e-07\\
2.1544e-06	9.80641071581688e-07\\
4.6416e-06	2.11276297715922e-06\\
1e-05	4.55178225901231e-06\\
};
\addlegendentry{$\en$,}		
\addplot [color=mycolor1, loosely dashed, line width=.5pt, mark=triangle*, mark options={solid, fill=mycolor1, mycolor1,rotate=90},mark size=4pt]
		table[row sep=crcr]{%
1e-14	1.23715599583333e-08\\
2.1544e-14	2.66532870131225e-08\\
4.6416e-14	5.74238258511125e-08\\
1e-13	1.23715582320237e-07\\
2.1544e-13	2.66532847835645e-07\\
4.6416e-13	5.74238238424928e-07\\
1e-12	1.23715580201148e-06\\
2.1544e-12	2.6653284564291e-06\\
4.6416e-12	5.74238096390502e-06\\
1e-11	1.23715531940061e-05\\
2.1544e-11	2.66532724026614e-05\\
4.6416e-11	5.74237956341811e-05\\
1e-10	0.000123715517918867\\
2.1544e-10	0.000266532709984063\\
4.6416e-10	0.000574237942151968\\
1e-09	0.0012371551643748\\
2.1544e-09	0.00266532708211884\\
4.6416e-09	0.00574237939024627\\
1e-08	0.0123715515496136\\
2.1544e-08	0.026653270435266\\
4.6416e-08	0.057423792162122\\
1e-07	0.123715507469232\\
2.1544e-07	0.26653266714726\\
4.6416e-07	0.57423774897374\\
1e-06	1.23715427338881\\
2.1544e-06	2.66532295232434\\
4.6416e-06	5.74236022641713\\
1e-05	12.3714626053369\\
		};
		\addlegendentry{$\ec$, }
		
\addplot [color=mycolor1, loosely dashed, line width=.5pt, mark=triangle, mark options={solid, fill=mycolor1, mycolor1,rotate=90},mark size=4pt]
		table[row sep=crcr]{%
1e-14	2.17686518847743e-07\\
2.1544e-14	4.68983836205616e-07\\
4.6416e-14	1.01041374588394e-06\\
1e-13	2.1768651884788e-06\\
2.1544e-13	4.68983836206252e-06\\
4.6416e-13	1.0104137458869e-05\\
1e-12	2.17686518849252e-05\\
2.1544e-12	4.68983836212618e-05\\
4.6416e-12	0.000101041374591644\\
1e-11	0.000217686518862966\\
2.1544e-11	0.000468983836276272\\
4.6416e-11	0.00101041374621191\\
1e-10	0.00217686519000108\\
2.1544e-10	0.0046898383691281\\
4.6416e-10	0.0101041374916657\\
1e-09	0.0217686520371534\\
2.1544e-09	0.0468983843278196\\
4.6416e-09	0.101041377871317\\
1e-08	0.217686534085785\\
2.1544e-08	0.468983906932065\\
4.6416e-08	1.01041407417922\\
1e-07	2.17686671228345\\
2.1544e-07	4.68984543471233\\
4.6416e-07	10.1041702884457\\
1e-06	21.7688042658802\\
2.1544e-06	46.8990908911401\\
4.6416e-06	101.044657598294\\
1e-05	217.701757448277\\
		};
        		\addlegendentry{upper bound, Example \ref{ex1}}
\addplot [color=mycolor2, dotted, line width=.5pt, mark=pentagon*, mark options={solid, fill=mycolor2, mycolor2,rotate=180},mark size=3.5pt]
  table[row sep=crcr]{%
1e-14	1.11278408774832e-14\\
2.1544e-14	2.15338388584813e-14\\
4.6416e-14	5.10952990443348e-14\\
1e-13	1.11184140951995e-13\\
2.1544e-13	2.41743830205111e-13\\
4.6416e-13	5.19118357894655e-13\\
1e-12	1.12224460410179e-12\\
2.1544e-12	2.41752041222094e-12\\
4.6416e-12	5.20955804095667e-12\\
1e-11	1.12172214637764e-11\\
2.1544e-11	2.41681224349732e-11\\
4.6416e-11	5.20725768967791e-11\\
1e-10	1.12182643383896e-10\\
2.1544e-10	2.41680687275213e-10\\
4.6416e-10	5.20698347187741e-10\\
1e-09	1.12180971277779e-09\\
2.1544e-09	2.41682578890366e-09\\
4.6416e-09	5.20699211965049e-09\\
1e-08	1.12180961790829e-08\\
2.1544e-08	2.4168270171079e-08\\
4.6416e-08	5.20699263066033e-08\\
1e-07	1.12180981269802e-07\\
2.1544e-07	2.41682660778326e-07\\
4.6416e-07	5.20698942682914e-07\\
1e-06	1.1218082100203e-06\\
2.1544e-06	2.41681917876301e-06\\
4.6416e-06	5.20695491503061e-06\\
1e-05	1.12179220159721e-05\\
};
\addlegendentry{$\en$,}		
		\addplot [color=mycolor2, dotted, line width=.5pt, mark=square*, mark options={solid, fill=mycolor2, mycolor2,rotate=45},mark size=3pt]
		table[row sep=crcr]{%
1e-14	1.53858943909226e-10\\
2.1544e-14	3.31486386227321e-10\\
4.6416e-14	7.14166071445348e-10\\
1e-13	1.53862006101212e-09\\
2.1544e-13	3.31479905397692e-09\\
4.6416e-13	7.14166079147461e-09\\
1e-12	1.53861938190586e-08\\
2.1544e-12	3.31480197157767e-08\\
4.6416e-12	7.14165667832648e-08\\
1e-11	1.53861980741283e-07\\
2.1544e-11	3.31480249287634e-07\\
4.6416e-11	7.14165758279192e-07\\
1e-10	1.53861979863364e-06\\
2.1544e-10	3.31480251771461e-06\\
4.6416e-10	7.14165769836119e-06\\
1e-09	1.5386198059432e-05\\
2.1544e-09	3.31480250481754e-05\\
4.6416e-09	7.1416576530034e-05\\
1e-08	0.000153861978529634\\
2.1544e-08	0.000331480240733309\\
4.6416e-08	0.00071416572008981\\
1e-07	0.0015386195753241\\
2.1544e-07	0.00331480143310308\\
4.6416e-07	0.00714165267905215\\
1e-06	0.0153861747651041\\
2.1544e-06	0.033147916915983\\
4.6416e-06	0.0714160746159471\\
1e-05	0.153859648879628\\
		};
		\addlegendentry{$\ec$, }
		
		\addplot [color=mycolor2, dotted, line width=.5pt, mark=square, mark options={solid, mycolor2,rotate=45},mark size=3pt]
		table[row sep=crcr]{%
1e-14	3.19979506152637e-09\\
2.1544e-14	6.89363848055297e-09\\
4.6416e-14	1.48521687575846e-08\\
1e-13	3.19979506152839e-08\\
2.1544e-13	6.89363848056233e-08\\
4.6416e-13	1.4852168757628e-07\\
1e-12	3.19979506154855e-07\\
2.1544e-12	6.8936384806559e-07\\
4.6416e-12	1.48521687580623e-06\\
1e-11	3.19979506175013e-06\\
2.1544e-11	6.89363848159155e-06\\
4.6416e-11	1.48521687624054e-05\\
1e-10	3.199795063766e-05\\
2.1544e-10	6.89363849094809e-05\\
4.6416e-10	0.000148521688058362\\
1e-09	0.000319979508392471\\
2.1544e-09	0.000689363858451352\\
4.6416e-09	0.00148521692401446\\
1e-08	0.00319979528551181\\
2.1544e-08	0.00689363952016784\\
4.6416e-08	0.0148521735832284\\
1e-07	0.0319979730138349\\
2.1544e-07	0.0689364887671825\\
4.6416e-07	0.148522170141383\\
1e-06	0.319981746017157\\
2.1544e-06	0.689374244293504\\
4.6416e-06	1.48526513303649\\
1e-05	3.20001905518014\\
		};
		\addlegendentry{upper bound, Example \ref{ex2}}
		
	\end{axis}
\end{tikzpicture}}
 	\caption{Example \ref{ex4}: Results for the minimal solution $\mm$ with $\alpha=0.499999999999999$. For this $\alpha$, Example \ref{ex1} yields $\mathrm{cond}(R_{\mm^\ast})=1.5382\cdot 10^{20}$, while Example \ref{ex2} yields $\mathrm{cond}(R_{\mm^\ast})=1.8159\cdot 10^{20}$.
    }
 	\label{figure6}
 \end{figure}
 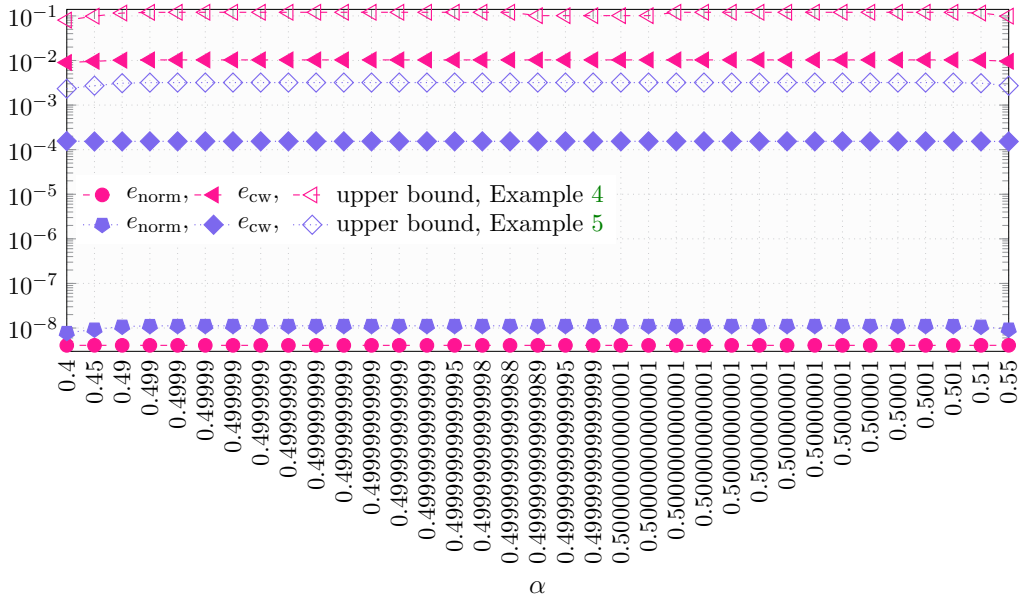
\begin{figure}[H]
 	\centerline{    
\definecolor{mycolor2}{RGB}{124,104,238}    
\definecolor{mycolor1}{RGB}{255,20,147}    
\begin{tikzpicture}[scale=.85]

	\begin{axis}[%
		width=5.8in,
		height=2.1in,
		xmin=2,
		xmax=36,
		ymode=log,
		ymin=3e-09,
		ymax=.14,
		scale only axis,
		xtick={1  ,2,  3,  4,  5,  6,  7,  8,  9, 10, 11, 12, 13, 14, 15, 16, 17, 18, 19, 20, 21, 22, 23, 24, 25, 26, 27, 28, 29, 30, 31, 32, 33, 34, 35, 36, 37, 38, 39, 40, 41},
		xticklabels={
			0.3,
			0.4,
			0.45,
			0.49,
			0.499,
			0.4999,
			0.49999,
			0.499999,
			0.4999999,
			0.49999999,
			0.499999999,
			0.4999999999,
			0.49999999999,
			0.499999999999,
			0.4999999999999,
			0.49999999999995,
			0.49999999999998,
			0.499999999999988,
			0.499999999999989,
			0.499999999999995,
			0.499999999999999,
			0.500000000000001,
			0.50000000000001,
			0.5000000000001,
			0.500000000001,
			0.50000000001,
			0.5000000001,
			0.500000001,
			0.50000001,
			0.5000001,
			0.500001,
			0.50001,
			0.5001,
			0.501,
			0.51,
			0.55,
			0.6,
			0.7,
			0.8,
			0.9,
			0.99,
		},
		xticklabel style={rotate=90},    
		yminorticks=true,
		xlabel={$\alpha$},
		xlabel style={font=\large},
		legend style={at={(0.009,0.312)}, anchor=south west, legend cell align=left,
			align=left,
			draw=none,
			font=\normalsize,
		},
		grid=major,
		grid style={dotted,gray!50},
		axis background/.style={fill=gray!2},
		legend columns=3
		]
		\addplot [color=mycolor1, loosely dashed, line width=.5pt, mark=*, mark options={solid, fill=mycolor1, mycolor1,rotate=90},mark size=3pt]
		table[row sep=crcr]{%
			1	4.01350407912851e-09\\
			2	4.08494125635081e-09\\
			3	4.09245083185765e-09\\
			4	4.08747409408698e-09\\
			5	4.08515252972851e-09\\
			6	4.08489696419477e-09\\
			7	4.08487147142077e-09\\
			8	4.08486868736869e-09\\
			9	4.08486798181138e-09\\
			10	4.08486926364339e-09\\
			11	4.08486755266842e-09\\
			12	4.08486796155392e-09\\
			13	4.08486855691345e-09\\
			14	4.08486851916663e-09\\
			15	4.0848696635796e-09\\
			16	4.08486926742368e-09\\
			17	4.08486912107811e-09\\
			18	4.08486762654699e-09\\
			19	4.08486823879724e-09\\
			20	4.08486746182059e-09\\
			21	4.08486720191526e-09\\
			22	4.08487021005662e-09\\
			23	4.0848862279382e-09\\
			24	4.08486836282775e-09\\
			25	4.08486776366303e-09\\
			26	4.08486867813087e-09\\
			27	4.08486794713476e-09\\
			28	4.08486720023966e-09\\
			29	4.08486836553848e-09\\
			30	4.08486848033542e-09\\
			31	4.08486871721313e-09\\
			32	4.08487138642909e-09\\
			33	4.08489685376565e-09\\
			34	4.08515236300449e-09\\
			35	4.0874741351244e-09\\
			36	4.09245124285909e-09\\
			37	4.08494125378802e-09\\
			38	4.01350401564053e-09\\
			39	3.84414138002043e-09\\
			40	3.54300798309586e-09\\
		};
		\addlegendentry{$\en$,}		
		\addplot [color=mycolor1, loosely dashed, line width=.5pt, mark=triangle*, mark options={solid, fill=mycolor1, mycolor1,rotate=90},mark size=4pt]
		table[row sep=crcr]{%
			1	0.00803763863252282\\
			2	0.00900280109193216\\
			3	0.0095997634058467\\
			4	0.0101548512006669\\
			5	0.0102912358552482\\
			6	0.0103051307213961\\
			7	0.0103065228253393\\
			8	0.0103066620619614\\
			9	0.0103066759858862\\
			10	0.0103066773782806\\
			11	0.0103066775175211\\
			12	0.0103066775314449\\
			13	0.0103066775328371\\
			14	0.0103066775329759\\
			15	0.0103066775329899\\
			16	0.0103066775329906\\
			17	0.010306677532991\\
			18	0.0103066775329919\\
			19	0.0103066775330103\\
			20	0.0103066775330004\\
			21	0.0103066775329938\\
			22	0.0103066775329893\\
			23	0.0103066775329689\\
			24	0.0103066775329901\\
			25	0.0103066775329762\\
			26	0.0103066775328369\\
			27	0.0103066775314445\\
			28	0.0103066775175207\\
			29	0.0103066773782811\\
			30	0.0103066759858855\\
			31	0.0103066620619613\\
			32	0.0103065228253396\\
			33	0.0103051307213962\\
			34	0.0102912358552482\\
			35	0.0101548512006667\\
			36	0.00959976340584663\\
			37	0.00900280109193232\\
			38	0.00803763863252322\\
			39	0.00726940804287745\\
			40	0.00661575732926487\\
		};
		\addlegendentry{$\ec$, }
		
		\addplot [color=mycolor1, loosely dashed, line width=.5pt, mark=triangle, mark options={solid, fill=mycolor1, mycolor1,rotate=90},mark size=4pt]
		table[row sep=crcr]{%
			1	0.0515050137178013\\
			2	0.0803550240689441\\
			3	0.0986443458894946\\
			4	0.115816067397194\\
			5	0.12005214007815\\
			6	0.120484020877179\\
			7	0.120527293341772\\
			8	0.120531621433748\\
			9	0.120532054251404\\
			10	0.120532097533254\\
			11	0.12053210186144\\
			12	0.120532102294259\\
			13	0.120532102337541\\
			14	0.120532102341869\\
			15	0.120532102342302\\
			16	0.120532102342324\\
			17	0.120532102342329\\
			18	0.120532102342349\\
			19	0.1013352098891\\
			20	0.1013352098891\\
			21	0.1013352098891\\
			22	0.1013352098891\\
			23	0.1013352098891\\
			24	0.120532102342301\\
			25	0.120532102341869\\
			26	0.120532102337541\\
			27	0.120532102294259\\
			28	0.12053210186144\\
			29	0.120532097533254\\
			30	0.120532054251404\\
			31	0.120531621433748\\
			32	0.120527293341811\\
			33	0.120484020877179\\
			34	0.12005214007815\\
			35	0.115816067397194\\
			36	0.0986443458894946\\
			37	0.0803550240689441\\
			38	0.0515050137178013\\
			39	0.0297922971165238\\
			40	0.0129801503301179\\
		};
		\addlegendentry{upper bound, Example \ref{ex1}}
		\addplot [color=mycolor2, dotted, line width=.5pt, mark=pentagon*, mark options={solid, fill=mycolor2, mycolor2,rotate=180},mark size=3.5pt]
		table[row sep=crcr]{%
			1	5.70295820286138e-09\\
			2	7.79995924608157e-09\\
			3	9.28780058977984e-09\\
			4	1.07886781204789e-08\\
			5	1.11740601652938e-08\\
			6	1.12136839639105e-08\\
			7	1.12176564522424e-08\\
			8	1.12180558233284e-08\\
			9	1.1218093953125e-08\\
			10	1.12180999186026e-08\\
			11	1.12180960844066e-08\\
			12	1.12180980799008e-08\\
			13	1.12180973066386e-08\\
			14	1.12180988919252e-08\\
			15	1.12180989266379e-08\\
			16	1.12180997195854e-08\\
			17	1.12180965341182e-08\\
			18	1.12180964452028e-08\\
			19	1.12180966206259e-08\\
			20	1.1218096355807e-08\\
			21	1.12180961792608e-08\\
			22	1.12180995399842e-08\\
			23	1.12181283706342e-08\\
			24	1.12180993162347e-08\\
			25	1.12180975700744e-08\\
			26	1.12180974420945e-08\\
			27	1.12180989038704e-08\\
			28	1.12180979626568e-08\\
			29	1.12180993767902e-08\\
			30	1.12180950278397e-08\\
			31	1.12180540094939e-08\\
			32	1.12176578099831e-08\\
			33	1.1213683586866e-08\\
			34	1.117405987184e-08\\
			35	1.07886784646596e-08\\
			36	9.28780087213981e-09\\
			37	7.79995924607267e-09\\
			38	5.70295791346428e-09\\
			39	4.36394937242333e-09\\
			40	3.51951194128637e-09\\
		};
		\addlegendentry{$\en$,}		
		\addplot [color=mycolor2, dotted, line width=.5pt, mark=square*, mark options={solid, fill=mycolor2, mycolor2,rotate=45},mark size=3pt]
		table[row sep=crcr]{%
			1	0.000178206433666146\\
			2	0.000155746916305981\\
			3	0.000152616831802425\\
			4	0.000153269569663297\\
			5	0.000153794936832426\\
			6	0.000153855196039718\\
			7	0.000153861299494122\\
			8	0.000153861910609354\\
			9	0.000153861971733369\\
			10	0.000153861977839552\\
			11	0.000153861978464869\\
			12	0.000153861978519797\\
			13	0.000153861978527512\\
			14	0.000153861978524052\\
			15	0.000153861978524013\\
			16	0.000153861978521584\\
			17	0.000153861978532088\\
			18	0.000153861978532371\\
			19	0.000153861978495835\\
			20	0.000153861978516078\\
			21	0.000153861978529633\\
			22	0.000153861978525838\\
			23	0.000153861978459145\\
			24	0.000153861978522601\\
			25	0.000153861978527897\\
			26	0.000153861978527082\\
			27	0.00015386197851622\\
			28	0.000153861978458633\\
			29	0.000153861977840996\\
			30	0.000153861971731011\\
			31	0.000153861910615649\\
			32	0.000153861299490592\\
			33	0.000153855196040644\\
			34	0.000153794936832828\\
			35	0.000153269569661825\\
			36	0.000152616831801689\\
			37	0.000155746916305409\\
			38	0.000178206433666953\\
			39	0.000240510560704092\\
			40	0.000449584216404427\\
		};
		\addlegendentry{$\ec$, }
		
		\addplot [color=mycolor2, dotted, line width=.5pt, mark=square, mark options={solid, mycolor2,rotate=45},mark size=3pt]
		table[row sep=crcr]{%
			1	0.00174785467925492\\
			2	0.00233320257596373\\
			3	0.00272120514927288\\
			4	0.00309542065934529\\
			5	0.00318905520382815\\
			6	0.00319862913505827\\
			7	0.00319987527139986\\
			8	0.00320001924506156\\
			9	0.00320003364288386\\
			10	0.00320003508267065\\
			11	0.00320003522664937\\
			12	0.00320003524104725\\
			13	0.00320003524248703\\
			14	0.00320003524263101\\
			15	0.00320003524264541\\
			16	0.00320003524264621\\
			17	0.00320003524264669\\
			18	0.00320003524264682\\
			19	0.00320003307821663\\
			20	0.00320003307821663\\
			21	0.00320003307821663\\
			22	0.00320003307821663\\
			23	0.00320003307821663\\
			24	0.00320003524264541\\
			25	0.00320003524263101\\
			26	0.00320003524248703\\
			27	0.00320003524104725\\
			28	0.00320003522664937\\
			29	0.00320003508267065\\
			30	0.00320003364288386\\
			31	0.00320001924506156\\
			32	0.0031998752714027\\
			33	0.00319862913505827\\
			34	0.00318905520382815\\
			35	0.00309542065934529\\
			36	0.00272120514927288\\
			37	0.00233320257596373\\
			38	0.00174785467925492\\
			39	0.00133333079297868\\
			40	0.00102925525908291\\
		};
		\addlegendentry{upper bound, Example \ref{ex2}}

	\end{axis}
\end{tikzpicture}%
}
 	\caption{Example \ref{ex4}: Results for the minimal solution $\mm$ with $\varepsilon=10^{-8}$.
    }
 	\label{figure8}
 \end{figure}
 	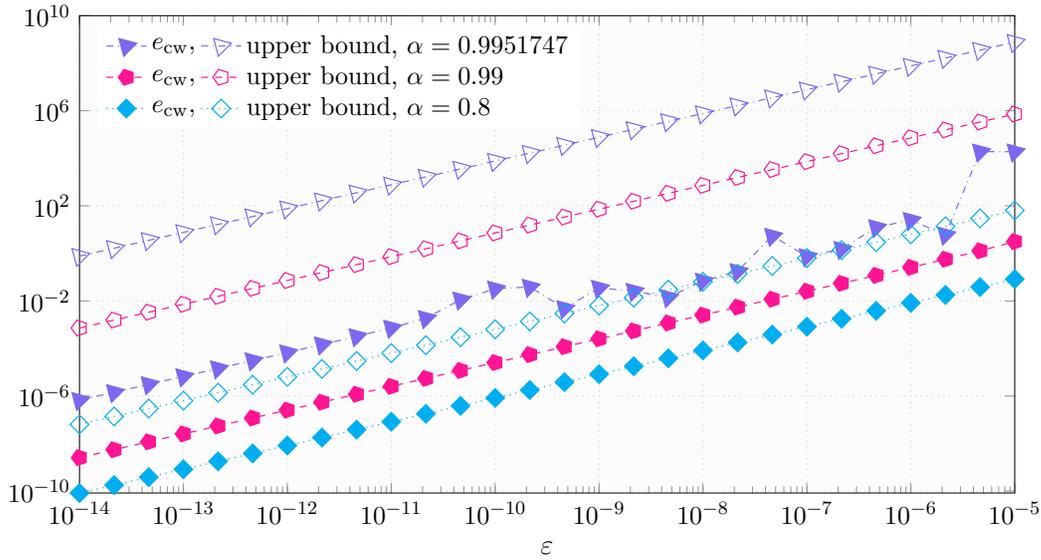
\begin{figure}[h]
 	 	\centerline{\definecolor{mycolor2}{RGB}{124,104,238}    
\definecolor{mycolor1}{RGB}{255,20,147} 
\begin{tikzpicture}[scale=.89]
	
	\begin{axis}[%
width=5.5in,
height=2.8in,
	scale only axis,
	xmode=log,
	xmin=1e-14,
	xmax=1e-05,
	xminorticks=true,
    xlabel={$\varepsilon$},
xlabel style={font=\large}, 
ymode=log,
	ymin=8.26107564057308e-11,
	ymax=10000000000,
	yminorticks=true,
		yminorticks=true,
		legend style={
			at={(0.02,0.98)},
			anchor=north west,
			legend cell align=left,
			align=left,
			draw=none,
			font=\normalsize,
		},
            grid=major,
 		grid style={dotted,gray!50},
 		axis background/.style={fill=gray!2},
		legend columns=2
		]
		\addplot [color=mycolor2, dashdotted, line width=.5pt, mark=triangle*, mark options={solid, fill=mycolor2, mycolor2,rotate=45},mark size=4.5pt]
table[row sep=crcr]{%
1e-14	6.65713085192857e-07\\
2.1544e-14	1.43432986412927e-06\\
4.6416e-14	3.0908382371738e-06\\
1e-13	6.66195779103371e-06\\
2.1544e-13	1.43664079092889e-05\\
4.6416e-13	3.1016792352837e-05\\
1e-12	6.71279818711261e-05\\
2.1544e-12	0.000146075907344217\\
4.6416e-12	0.000321932032821603\\
1e-11	0.000732818152703823\\
2.1544e-11	0.00187215404145333\\
4.6416e-11	0.011268898841526\\
1e-10	0.0324224535113928\\
2.1544e-10	0.0387419871881065\\
4.6416e-10	0.00488532096714623\\
1e-09	0.0349209377615315\\
2.1544e-09	0.0260429435074099\\
4.6416e-09	0.0141929073112445\\
1e-08	0.069605607142109\\
2.1544e-08	0.177130024356941\\
4.6416e-08	5.40192034283143\\
1e-07	0.750866059319094\\
2.1544e-07	1.51625831600001\\
4.6416e-07	12.5309177220278\\
1e-06	25.0733015285366\\
2.1544e-06	6.06106674436685\\
4.6416e-06	19017.0794588182\\
1e-05	19017.7353250032\\
};
\addlegendentry{$\ec$,}

\addplot [color=mycolor2, dashdotted, line width=.5pt, mark=triangle, mark options={solid, mycolor2,rotate=45},mark size=4.5pt]
table[row sep=crcr]{%
1e-14	0.766482694816301\\
2.1544e-14	1.65131031771237\\
4.6416e-14	3.55770607626025\\
1e-13	7.66482694816784\\
2.1544e-13	16.5131031771461\\
4.6416e-13	35.5770607627065\\
1e-12	76.6482694821613\\
2.1544e-12	165.131031773703\\
4.6416e-12	355.770607637469\\
1e-11	766.482694869901\\
2.1544e-11	1651.31031796115\\
4.6416e-11	3557.70607741503\\
1e-10	7664.82695352785\\
2.1544e-10	16513.1032020243\\
4.6416e-10	35577.0608781851\\
1e-09	76648.2700181626\\
2.1544e-09	165131.03426152\\
4.6416e-09	355770.619185323\\
1e-08	766482.748470038\\
2.1544e-08	1651310.56674295\\
4.6416e-08	3557707.23220068\\
1e-07	7664832.31354325\\
2.1544e-07	16513128.080221\\
4.6416e-07	35577176.3569222\\
1e-06	76648806.0214273\\
2.1544e-06	165133522.098437\\
4.6416e-06	355782167.2315\\
1e-05	766536350.521139\\
};
\addlegendentry{upper bound, $\alpha=0.9951747$}

		\addplot [color=mycolor1, dashed, line width=.5pt, mark=pentagon*, mark options={solid, fill=mycolor1, mycolor1,rotate=90},mark size=3.5pt]
		table[row sep=crcr]{%
1e-14	2.54267319722923e-09\\
2.1544e-14	5.47784551490036e-09\\
4.6416e-14	1.18018384499989e-08\\
1e-13	2.54262720203625e-08\\
2.1544e-13	5.47783378422448e-08\\
4.6416e-13	1.18018521247691e-07\\
1e-12	2.54262546368783e-07\\
2.1544e-12	5.47783208450573e-07\\
4.6416e-12	1.18018509517016e-06\\
1e-11	2.54262576898966e-06\\
2.1544e-11	5.47783396306524e-06\\
4.6416e-11	1.1801858522858e-05\\
1e-10	2.54262932818566e-05\\
2.1544e-10	5.47785042754544e-05\\
4.6416e-10	0.000118019349331302\\
1e-09	0.000254266481571037\\
2.1544e-09	0.000547801515411225\\
4.6416e-09	0.00118026996147639\\
1e-08	0.00254301981222874\\
2.1544e-08	0.0054796634981291\\
4.6416e-08	0.0118103573376568\\
1e-07	0.0254658071913567\\
2.1544e-07	0.0549625690244192\\
4.6416e-07	0.118880462673226\\
1e-06	0.2583328285185\\
2.1544e-06	0.567420055907224\\
4.6416e-06	1.28012907689453\\
1e-05	3.14627338574067\\
		};
		\addlegendentry{$\ec$,}
		
		\addplot [color=mycolor1, dashed, line width=.5pt, mark=pentagon, mark options={solid, fill=mycolor1, mycolor1,rotate=90},mark size=3.5pt]
		table[row sep=crcr]{%
1e-14	0.00072418884891547\\
2.1544e-14	0.00156019245610361\\
4.6416e-14	0.0033613949611269\\
1e-13	0.00724188848915926\\
2.1544e-13	0.0156019245610573\\
4.6416e-13	0.0336139496113673\\
1e-12	0.0724188848920488\\
2.1544e-12	0.156019245612691\\
4.6416e-12	0.336139496123502\\
1e-11	0.724188848966112\\
2.1544e-11	1.56019245633867\\
4.6416e-11	3.36139496221797\\
1e-10	7.24188849422351\\
2.1544e-10	15.6019245845627\\
4.6416e-10	33.6139497204738\\
1e-09	72.4188853984741\\
2.1544e-09	156.019247963233\\
4.6416e-09	336.139507034157\\
1e-08	724.18889960864\\
2.1544e-08	1560.1926913929\\
4.6416e-08	3361.39605328355\\
1e-07	7241.89355847794\\
2.1544e-07	15601.9480900023\\
4.6416e-07	33614.0588271951\\
1e-06	72419.3918255466\\
2.1544e-06	156021.598523481\\
4.6416e-06	336150.417869228\\
1e-05	724239.543945361\\
		};
		\addlegendentry{upper bound, $\alpha=0.99$}
		\addplot [color=cyan, dotted, line width=.5pt, mark=square*, mark options={solid, fill=cyan, cyan,rotate=45},mark size=2.8pt]
		table[row sep=crcr]{%
1e-14	8.26107564057308e-11\\
2.1544e-14	1.77975204494694e-10\\
4.6416e-14	3.83443632100817e-10\\
1e-13	8.26100971551182e-10\\
2.1544e-13	1.77975070524266e-09\\
4.6416e-13	3.83443091461859e-09\\
1e-12	8.26100990081453e-09\\
2.1544e-12	1.77975184211139e-08\\
4.6416e-12	3.83443023842346e-08\\
1e-11	8.26100966517697e-08\\
2.1544e-11	1.77975191174251e-07\\
4.6416e-11	3.83443024705429e-07\\
1e-10	8.26100964839102e-07\\
2.1544e-10	1.77975191921395e-06\\
4.6416e-10	3.83443023334509e-06\\
1e-09	8.26100961807618e-06\\
2.1544e-09	1.77975190266074e-05\\
4.6416e-09	3.83443015946856e-05\\
1e-08	8.26100927296869e-05\\
2.1544e-08	0.000177975174245214\\
4.6416e-08	0.000383442941632095\\
1e-07	0.000826100582358604\\
2.1544e-07	0.0017797501414439\\
4.6416e-07	0.00383442198483957\\
1e-06	0.0082609713300325\\
2.1544e-06	0.0177973413158975\\
4.6416e-06	0.0383434767214706\\
1e-05	0.0826062641481303\\
		};
		\addlegendentry{$\ec$,}
		
	\addplot [color=cyan, dotted, line width=.5pt, mark=square, mark options={solid, fill=cyan, cyan,rotate=45},mark size=2.8pt]
		table[row sep=crcr]{%
1e-14	6.39465805852345e-08\\
2.1544e-14	1.3776651321284e-07\\
4.6416e-14	2.968144484445e-07\\
1e-13	6.39465805852748e-07\\
2.1544e-13	1.37766513213027e-06\\
4.6416e-13	2.96814448445368e-06\\
1e-12	6.39465805856777e-06\\
2.1544e-12	1.37766513214897e-05\\
4.6416e-12	2.96814448454048e-05\\
1e-11	6.39465805897063e-05\\
2.1544e-11	0.000137766513233596\\
4.6416e-11	0.000296814448540842\\
1e-10	0.000639465806299926\\
2.1544e-10	0.00137766513420583\\
4.6416e-10	0.00296814449408789\\
1e-09	0.00639465810328561\\
2.1544e-09	0.0137766515290449\\
4.6416e-09	0.0296814458088261\\
1e-08	0.0639465850614908\\
2.1544e-08	0.137766533989113\\
4.6416e-08	0.296814544882996\\
1e-07	0.639466253478524\\
2.1544e-07	1.37766720975903\\
4.6416e-07	2.96815412831763\\
1e-06	6.39470282128928\\
2.1544e-06	13.7768590858038\\
4.6416e-06	29.682409246188\\
1e-05	63.9510570057415\\
		};
		\addlegendentry{upper bound, $\alpha=0.8$}
		
	\end{axis}
\end{tikzpicture}}
 	\caption{Example \ref{ex4}: Results for the stochastic solution $\ss$, where $\PP$ and $\vv$ are taken from Example \ref{ex2}.}
 	\label{figure7}
 \end{figure}
 \end{example}
It is interesting to note that the theoretical perturbation bound that we obtained are quite loose, and that in practice they seem to hold also for stochastic solutions, even though the proof only works for minimal solutions.
\section{Conclusions}\label{sec:conclusions}
In this work, we show that the additional constraint $\mathbf{1}_n^T\aa = 1-\alpha$, $\mathbf{1}_n^T \BB_{(1)} = \alpha \mathbf{1}_{n^2}^T$ in~\eqref{mlpr} allows one to compute solutions to this equation in a more stable way than perturbation bounds based on the norm of the Jacobian $\norm{R_{\xx}}$ would predict.

Together with these theoretical results, we have obtained an interesting new definition of partial inverse of an M-matrix, which stays bounded even when the matrix approaches singularity, and new interesting algorithms for~\eqref{eq01} and~\eqref{mlpr}. To our knowledge, also the subtraction-free formulation of the Newton method with the formula $\mathbf{r} = \BB \hh^2$ in Algorithm~\ref{algo1} is novel, in addition to the block Jacobi algorithm and its variant.

It would be interesting to extend the scope of our research and show if similar results hold also for other problems which have constraints on stochasticity and sum of entries, for instance the one considered in~\cite{xue2012accurate}, which is a special case of~\eqref{eq01} appearing in queuing theory.


%

\bibliographystyle{siamplain}
\bibliography{references}

\begin{thebibliography}{10}

\bibitem{Abdesselam}
{\sc A.~Abdesselam}, {\em The {G}rassmann--{B}erezin calculus and theorems of
  the matrix-tree type}, Adv. Appl. Math., 33 (2004), pp.~51--70,
  \url{https://doi.org/10.1016/j.aam.2003.07.002}.

\bibitem{alfa}
{\sc A.~S. Alfa, J.~Xue, and Q.~Ye}, {\em Entrywise perturbation theory for
  diagonally dominant {M}-matrices with applications}, Numer. Math., 90 (2002),
  pp.~401--414, \url{https://doi.org/10.1007/s002110100289}.

\bibitem{batage}
{\sc V.~Batagelj and A.~Mrvar}, {\em Pajek datasets}, 2006,
  \url{http://vlado.fmf.uni-lj.si/pub/networks/data}.

\bibitem{benson}
{\sc A.~R. Benson, D.~F. Gleich, and J.~Leskovec}, {\em Tensor spectral
  clustering for partitioning higher-order network structures}, in Proceedings
  of the 2015 SIAM International Conference on Data Mining, SIAM, 2015,
  pp.~118--126, \url{https://doi.org/0.1137/1.9781611974010.14}.

\bibitem{benson2}
{\sc A.~R. Benson, D.~F. Gleich, and L.-H. Lim}, {\em The spacey random walk: a
  stochastic process for higher-order data}, SIAM Rev., 59 (2017),
  pp.~321--345, \url{https://doi.org/10.1137/16M1074023}.

\bibitem{chaiken}
{\sc S.~Chaiken}, {\em A combinatorial proof of the all minors matrix tree
  theorem}, SIAM J. Alg. Disc. Meth., 3 (1982), pp.~319--329,
  \url{https://doi.org/10.1137/0603033}.

\bibitem{cipolla2020extrapolation}
{\sc S.~Cipolla, M.~Redivo-Zaglia, and F.~Tudisco}, {\em Extrapolation methods
  for fixed-point multilinear {P}age{R}ank computations}, Numer. Linear Algebra
  Appl., 27 (2020), p.~e2280, \url{https://doi.org/10.1002/nla.2280}.

\bibitem{Davis}
{\sc T.~A. Davis and Y.~Hu}, {\em The university of {F}lorida sparse matrix
  collection}, ACM Transactions on Mathematical Software (TOMS), 38 (2011),
  pp.~1--25, \url{https://doi.org/10.1145/2049662.2049663}.

\bibitem{gleich2015}
{\sc D.~F. Gleich, L.-H. Lim, and Y.~Yu}, {\em Multilinear {P}age{R}ank}, SIAM
  J. Matrix Anal. Appl., 36 (2015), pp.~1507--1541,
  \url{https://doi.org/10.1137/140985160}.

\bibitem{grindrod}
{\sc P.~Grindrod and T.~Lee}, {\em Comparison of social structures within
  cities of very different sizes}, Royal Society Open Science, 3 (2016),
  p.~150526, \url{https://doi.org/10.1098/rsos.150526}.

\bibitem{higham-accuracy}
{\sc N.~J. Higham}, {\em Accuracy and stability of numerical algorithms.},
  Philadelphia, PA: SIAM, 2nd ed.~ed., 2002,
  \url{https://doi.org/10.1137/1.9780898718027}.

\bibitem{lai2023anderson}
{\sc F.~Lai, W.~Li, X.~Peng, and Y.~Chen}, {\em Anderson accelerated
  fixed-point iteration for multilinear {P}age{R}ank}, Numer. Linear Algebra
  Appl., 30 (2023), p.~e2499, \url{https://doi.org/10.1002/nla.2499}.

\bibitem{meini2018perron}
{\sc B.~Meini and F.~Poloni}, {\em Perron-based algorithms for the multilinear
  {P}age{R}ank}, Numer. Linear Algebra Appl., 25 (2018), p.~e2177,
  \url{https://doi.org/10.1002/nla.2177}.

\bibitem{moler}
{\sc C.~B. Moler}, {\em Numerical computing with \textsc{Matlab}}, SIAM, 2004.

\bibitem{Ocinneide}
{\sc C.~A. O'Cinneide}, {\em Relative-error bounds for the {LU} decomposition
  via the {GTH} algorithm}, Numer. Math., 73 (1996), pp.~507--519,
  \url{https://doi.org/10.1007/s002110050203}.

\bibitem{qve}
{\sc F.~Poloni}, {\em Quadratic vector equations}, Linear Algebra Appl., 438
  (2013), pp.~1627--1644, \url{https://doi.org/10.1016/j.laa.2011.05.036}.

\bibitem{saad}
{\sc Y.~Saad}, {\em Iterative {M}ethods for {S}parse {L}inear {S}ystems}, SIAM,
  2003.

\bibitem{stewartsun}
{\sc G.~W. Stewart and J.-g. Sun}, {\em Matrix perturbation theory}, Boston
  etc.: Academic Press, Inc., 1990.

\bibitem{stoer}
{\sc J.~Stoer, R.~Bulirsch, R.~Bartels, W.~Gautschi, and C.~Witzgall}, {\em
  Introduction to {N}umerical {A}nalysis}, vol.~1993, Springer, 1980.

\bibitem{tisseur}
{\sc F.~Tisseur}, {\em Newton's method in floating point arithmetic and
  iterative refinement of generalized eigenvalue problems}, SIAM J. Matrix
  Anal. Appl., 22 (2001), pp.~1038--1057,
  \url{https://doi.org/10.1137/S0895479899359837},
  \url{eprints.maths.manchester.ac.uk/979/1/tisseur3.pdf}.

\bibitem{xxl2}
{\sc J.~Xue, S.~Xu, and R.-C. Li}, {\em Accurate solutions of {M}-matrix
  algebraic {Riccati} equations}, Numer. Math., 120 (2012), pp.~671--700,
  \url{https://doi.org/10.1007/s00211-011-0421-0}.

\bibitem{xue2012accurate}
{\sc J.~Xue, S.~Xu, and R.-C. Li}, {\em Accurate solutions of {M}-matrix
  {S}ylvester equations}, Numer. Math., 120 (2012), pp.~639--670,
  \url{https://doi.org/10.1007/s00211-011-0420-1}.

\end{thebibliography}

\end{document}


\maketitle

\section{A detailed example}

Here we include some equations and theorem-like environments to show
how these are labeled in a supplement and can be referenced from the
main text.
Consider the following equation:
\begin{equation}
  \label{eq:suppa}
  a^2 + b^2 = c^2.
\end{equation}
You can also reference equations such as \cref{eq:matrices,eq:bb} 
from the main article in this supplement.

\lipsum[100-101]

\begin{theorem}
An example theorem.
\end{theorem}

\lipsum[102]
 
\begin{lemma}
An example lemma.
\end{lemma}

\lipsum[103-105]

Here is an example citation: \cite{KoMa14}.

\section[Proof of Thm]{Proof of \cref{thm:bigthm}}
\label{sec:proof}

\lipsum[106-112]

\section{Additional experimental results}
\Cref{tab:foo} shows additional
supporting evidence. 

\begin{table}[htbp]
\footnotesize
  \caption{Example table.}  \label{tab:smfoo}
\begin{center}
  \begin{tabular}{|c|c|c|} \hline
   Species & \bf Mean & \bf Std.~Dev. \\ \hline
    1 & 3.4 & 1.2 \\
    2 & 5.4 & 0.6 \\ \hline
  \end{tabular}
\end{center}
\end{table}

\bibliographystyle{siamplain}
\bibliography{references}